\documentclass[a4paper,fleqn]{cas-sc}

\usepackage[numbers]{natbib}

\usepackage{cite}
\usepackage{algorithm}
\usepackage{algorithmic}
\usepackage{subfigure}
\usepackage{appendix}
\usepackage{graphicx}
\usepackage{caption}

\def\tsc#1{\csdef{#1}{\textsc{\lowercase{#1}}\xspace}}
\tsc{WGM}
\tsc{QE}
\tsc{EP}
\tsc{PMS}
\tsc{BEC}
\tsc{DE}

\begin{document}
\let\WriteBookmarks\relax
\def\floatpagepagefraction{1}
\def\textpagefraction{.001}

\shorttitle{Few-shot Learning for Inference in Medical Imaging}

\shortauthors{J. Liu et~al.}

\title [mode = title]{Few-shot Learning for Inference in Medical Imaging with Subspace Feature Representations}                      


\author{Jiahui Liu}[style=chinese]\fnref{equal}
\ead{jl4f19@soton.ac.uk}
\cormark[1]
\author{Keqiang Fan}[style=chinese]\fnref{equal}
\ead{k.fan@soton.ac.uk}
\cormark[1]
\cortext[cor1]{Corresponding author}
\fntext[equal]{Equal contribution}
\author{Xiaohao Cai}[style=chinese]
\ead{x.cai@soton.ac.uk}
\author{Mahesan Niranjan}
\ead{mn@ecs.soton.ac.uk}

\affiliation{organization={School of Electronics and Computer Science,
         University of Southampton}, 
    city={Southampton},
    postcode={SO17 1BJ}, 
    country={UK}}




\begin{abstract}
Unlike the field of visual scene recognition where tremendous advances have taken place due to the availability of very large datasets to train deep neural networks, inference from medical images is often hampered by the fact that only small amounts of data may be available. When working with very small dataset problems, of the order of a few hundred items of data, the power of deep learning may still be exploited by using a model pre-trained on natural images as a feature extractor and carrying out classic pattern recognition techniques in this feature space, the so-called few-shot learning problem. In regimes where the dimension of this feature space is comparable to or even larger than the number of items of data, dimensionality reduction is a necessity and is often achieved by principal component analysis, i.e., singular value decomposition (SVD). In this paper, noting the inappropriateness of using SVD for this setting, we usher in and explore two alternatives
based on discriminant analysis and non-negative matrix factorization (NMF). Using $14$ different datasets spanning $11$ distinct disease types, we demonstrate that discriminant subspaces at low dimensions achieve significant improvements over SVD-based subspaces and the original feature space. We also show that NMF at modest dimensions is a competitive alternative to SVD in this setting.
\end{abstract}


\begin{keywords}
Medical imaging \sep Few-shot learning \sep Classification \sep Discriminant analysis \sep  PCA \sep Non-negative matrix factorization \sep Dimensionality reduction
\end{keywords}

\maketitle

\section{Introduction}\label{sec:introduction}
Impressive empirical performances have been reported in the field of computer vision in recent years, starting from a step improvement reported in the ImageNet challenge \citep{krizhevsky2012imagenet}. This and subsequent work has used very large neural network architectures, notably their depth, with parameter estimation carried out using equally large datasets. It is common in current computer vision literature to train models with tens of millions of parameters and use datasets of similar sizes. Much algorithmic development to control the complexity of such massive models and to incorporate techniques to handle systematic variability has been developed. 
Our curiosity about mammalian vision \citep{fukushima1982neocognitron,pinto2008real} and commercial applications such as self-driving cars and robot navigation \citep{GUPTA2021100057,zhu2017target} has driven the computer vision field. The interest in automatic diagnosis has reached a level of comparing artificial intelligence-based methods against human clinicians \citep{liu2019comparison,pesapane2018artificial}.
However, compared with natural images, the application of deep learning in the medical domain poses more challenges, such as causality \citep{castro2020causality}, uncertainty \citep{kendall2017uncertainties}, and the need to integrate clinical information along with features extracted from images \citep{LUNDERVOLD2019102}. A particular issue with image-based inference in the medical field is data availability \citep{kim2017few}. Often, the number of images available in the medical domain is orders of magnitude smaller than what is state-of-the-art in computer vision. Compared with other domains, due to privacy concerns and the prevalence of adverse medical conditions, most of the medical datasets only contain thousands or even hundreds of images, such as brain imaging \citep{althnian2021impact}.

The focus of this paper is on data sparsity/scarcity. Naturally, if we had access to hundreds of thousands of labelled medical images, as might be the case with X-rays and optometry, training a deep neural network from scratch using all the recent methodological advances is the way forward. When the number of images is in the thousands, the strategy of transfer learning is suitable for the medical data by fine-tuning the weights generated from pre-trained natural images. While the scheme is appealing, available empirical evidence for transfer learning is contradictory in the medical field. For example, on a chest X-rays problem, Raghu {\em et al.} \citep{raghu2019transfusion} found no significant improvement with the popular ResNet trained on ImageNet as source architecture; more positive results are reported for endoscopy image recognition  \citep{CAROPPO2021101852}. Another example may be the weakly supervised learning methods \citep{zhang2021weakly}, whose performance is yet to be seen in medical diagnosis.

Our interest is in a regime of even smaller amounts of data than is needed to fine-tune a pre-trained model with transfer learning. 
This regime is referred to as ``few-shot learning'' \citep{huang2021unsupervised,peng2019few,tang2020blockmix,wang2020generalizing}, and is appropriate for dataset sizes of the order of a few hundred or even down to a few tens \citep{argueso2020few,quellec2020automatic}. 
Few-shot learning works can be divided into different categories -- data, model and algorithm \citep{wang2020generalizing}. Most contemporary few-shot learning techniques rely on models and algorithms with fine-tuned parameters based on available data \citep{8935497,9311786}. Data augmentation technology and manifold space \citep{das2020few} have also drawn some attention.
Unlike these methods, we in this paper explore few-shot problems from the traditional machine learning perspective by using a pre-trained deep neural network as a feature extractor. In detail, each image is mapped into a fixed dimensional feature space, the dimensions of which, say $M$, are defined by the number of neurons in the penultimate fully connected layer of the network, typically $512$ or $1024$ for the popular architectures. Then we are in a regime where the number of items of data, say $N$, is comparable to or even smaller than the dimension of the feature space (i.e., the $N < M$ problem in statistical inference language \citep{tibshirani1996regression}), necessitating techniques for dimensionality reduction.

Subspace methods for reducing the dimensionality of data have a long and rich history. They fall under the group of methods known as structured low-rank approximation methods \citep{markovsky2012low, shetta2021convex,markovsky2010approximate}. The basic intuition is a data matrix, $\boldsymbol{Y} \in \mathbb{R}^{N \times M}$, consisting of $N$ items of data in $M$ dimensional features, is usually not full rank. This is due to correlations along either of the axes. In the medical context, profiles of patients (i.e. data) may show strong similarities. Along the features axis, some features that have been gathered may be derivable from others. In these situations, we can find low-rank approximations by factorising $\boldsymbol{Y}$, and additionally impose structural constraints on the factors either from prior knowledge or for mathematical convenience. Popular approaches like principal component analysis (PCA) \citep{papailiopoulos2013sparse} and non-negative matrix factorization (NMF) \citep{lee1999learning,li2017robust} impose orthogonality and non-negativity constraints on the factors, respectively. Returning to few-shot learning with pre-trained deep neural network as feature extractors and encountering $N < M$ problems, pattern recognition problems are known to suffer the ``curse of dimensionality''. Hence dimensionality reduction techniques are required. The most popular technique used hitherto in the literature is PCA, implemented via singular value decomposition (SVD) \citep{quellec2020automatic,raghu2017svcca,wu2021generalized,9550610}. Despite its popularity, PCA has a fundamental weakness in that it is a variance-preserving low-rank approximation technique, more suitable for data that is uni-modal and Gaussian distributed. In the case of classification problems, however, the feature space is necessarily multi-modal with at least as many modes as the number of classes in the problem.

The basic premise of this work is the need for dimensionality reduction in the feature space and that SVD ignores multi-modal data structure. We, for the first time, usher in and explore two alternatives -- {\em discriminant analysis (DA) subspace} and {\em NMF subspace} -- to SVD in few-shot learning on medical imaging with multi-modal structure in the data.  The DA subspace introduces the well-known Fisher linear discriminant analysis (FDA) and its multi-dimensional extensions \citep{foley1975optimal}. The NMF  \citep{lee1999learning} and the supervised NMF (SNMF) \citep{leuschner2019supervised} (where class label information can be injected into the factorization loss function) subspaces focus on the part-based representation with sparsity.
A detailed comparison between these subspace representations, including feature selection techniques \citep{li2015robust}, is conducted.
Validating on $14$ datasets spanning $11$ medical classification tasks with four distinct imaging modalities, we achieve statistically significant improvements in classification accuracy in the subspaces compared to the original high-dimensional feature space, with persuasive results on DA and NMF subspaces as viable alternatives to SVD.

The remainder of this paper is organized as follows. In the next section, we mainly recall the subspace representation methods, i.e., SVD, DA and NMF subspaces. The few-shot learning methodology/scheme on subspace feature representations including the experimental settings in sufficient detail to facilitate reproduction of the work is provided in Section \ref{Sec:fl-art}. In Section \ref{sec:data}, we give succinct descriptions of the datasets used. Section \ref{sec:results} presents the key results of the experimental work. A further discussion is conducted in Section \ref{Sec:dis}, followed by conclusion in Section \ref{Sec:con}. Some additional details regarding method derivations and extra results are provided in Appendix.

\section{Subspace Representation}\label{Sec:subspace-rep}
\subsection{Basic Notations}
Given $N$ samples $\boldsymbol{y}_i = (y_{i1}, y_{i2}, \cdots, y_{iM})^\top \in \mathbb{R}^M, 1\le i \le N$,
we form a data matrix $\boldsymbol{Y} = (\boldsymbol{y}_1, \boldsymbol{y}_2, \cdots, \boldsymbol{y}_N)^\top \in \mathbb{R}^{N\times M}$, where $M$ is the number of features of every sample. Suppose that these $N$ samples belong to $C$ different classes, namely $\boldsymbol{\Lambda}_j$, and their cardinality $|\boldsymbol{\Lambda}_j| = N_j$, $1\le j \le C$. Let $\boldsymbol{y}^j_k$ represent the $k$-th sample in class $\boldsymbol{\Lambda}_j$. Clearly, $N = \sum_{j=1}^C N_j$, $\boldsymbol{\Lambda}_j = \{\boldsymbol{y}^j_k\}_{k=1}^{N_j}$ and $\{\boldsymbol{y}_i\}_{i=1}^N = \bigcup_{j=1}^C \{\boldsymbol{y}^j_k\}_{k=1}^{N_j}$. Let $\bar{\boldsymbol{y}}$ and $\bar{\boldsymbol{y}}_j$ respectively be the mean of the whole samples and the samples in class $j$, i.e., $\bar{\boldsymbol{y}} = \frac{1}{N} \sum_{i=1}^N \boldsymbol{y}_i$, $\bar{\boldsymbol{y}}_j = \frac{1}{N_j} \sum_{\boldsymbol{y} \in \boldsymbol{\Lambda}_j} \boldsymbol{y}$, $1\le j \le C$.

Let $\boldsymbol{S}^j_{\rm W}$ represent the intra-class scatter for class $j$, i.e.,
\begin{equation}
\boldsymbol{S}^j_{\rm W}= \sum_{k=1}^{N_{j}}(\boldsymbol{y}_k^j-\bar{\boldsymbol{y}}_j) (\boldsymbol{y}_k^j-\bar{\boldsymbol{y}}_j)^{\top}, \ \ 1\le j \le C.
\end{equation}
Then the inter- and intra-class scatters, denoted as $\boldsymbol{S}_{\rm B}$ and $\boldsymbol{S}_{\rm W}$, respectively, read
\begin{equation}
\boldsymbol{S}_{\rm B}  =\sum_{j=1}^C(\bar{\boldsymbol{y}}_j - \bar{\boldsymbol{y}})(\bar{\boldsymbol{y}}_j - \bar{\boldsymbol{y}})^{\top}, \quad
\boldsymbol{S}_{\rm W}  =\sum_{j=1}^C \boldsymbol{S}^j_{\rm W}.
\end{equation}
Specifically, for the binary case, i.e., $C=2$, we also name
$\tilde{\boldsymbol{S}}_{\rm B}$ and $\tilde{\boldsymbol{S}}_{\rm W}$ as the inter- and intra-class scatters, i.e.,
\begin{equation} \label{eqn:bi-intra-inter}
\tilde{\boldsymbol{S}}_{\rm B} = \boldsymbol{s}_{\rm b} \boldsymbol{s}_{\rm b}^\top, \quad  
\tilde{\boldsymbol{S}}_{\rm W} = \beta \boldsymbol{S}_{\rm W}^{1}+(1-\beta) \boldsymbol{S}_{\rm W}^{2},
\end{equation}
where $\boldsymbol{s}_{\rm b} = \bar{\boldsymbol{y}}_1 - \bar{\boldsymbol{y}}_2$ and $\beta=(N_{2}-1)/(N_{1}+N_{2}-2)$.

\subsection{Feature Selection}

Feature selection has been used in medical imaging \citep{tang2020cart}. It is the process of extracting a subset of relevant features by eliminating redundant or unnecessary information for model development \citep{li2015robust}. There are several types of feature selection techniques, including supervised \citep{huang2015supervised,kursa2010feature}, semi-supervised \citep{li2021semi}, and unsupervised methods \citep{li2015unsupervised}.  
For example, the Boruta algorithm \citep{kursa2010feature}, one of the supervised feature selection methods, selects features by shuffling features of the data and calculating the feature correlations based on classification loss. The approach has also been used to classify medical images \citep{tang2020cart}.

\subsection{Singular Value Decomposition}
SVD is the most common type of matrix decomposition, which can decompose either a square or rectangle matrix. The SVD of the matrix $\boldsymbol{Y}$ can be represented as 
$\boldsymbol{Y}=\boldsymbol{U} \boldsymbol{\Sigma} \boldsymbol{V}^{\top}$,
where $\boldsymbol{U} \in \mathbb{R}^{N \times N}$ and $\boldsymbol{V} \in \mathbb{R}^{M \times M}$ are orthogonal matrices, and $\boldsymbol{\Sigma} \in \mathbb{R}^{N \times M}$ is a diagonal matrix whose diagonal consists of singular values.
The singular values are generally ordered and it is well known that in most real-world problems they reduce quickly to zero, i.e., typically the first 10\% or even 1\% of the largest singular values could account for more than 99\% of the sum of all the singular values. Therefore, the singular vectors corresponding to the top $p \ll \min\{M, N\}$ largest singular values compose the transformation matrix for the most representative subspace. Meanwhile, the variance preserving property of SVD is extremely effective in data compression and widely employed in deep learning tasks, especially when the data is uni-modal. For example, SVD has been used to compress features taken at different layers to compare learning dynamics across layers as well as across models \citep{raghu2017svcca}. 

\subsection{Discriminant Subspaces} \label{Sec:sub-DA}
It is usually possible to design logic based on the statistics of a design set that achieves a very high recognition rate if the original set of features is well chosen. Discriminant vectors for DA can reduce the error rate and solve the discrimination portion of the task \citep{foley1975optimal}. 
Since the discriminant vector transformation aims to reduce dimensionality while retaining discriminatory information, sophisticated pattern recognition techniques that were either computationally impractical or statistically insignificant in the original high-dimensional space could be possible in the new and low-dimensional space.
The intuitive assumption is that features based on discrimination are better than that based on fitting or describing the data. In what follows, we present different approaches of obtaining discriminant vectors for multiclass and binary classification problems.

\subsubsection{{Multiclass classification Problem}}
The aim of the multiclass DA is to discover a linear transformation which lowers the dimensionality of an $M$-dimensional statistical model with $C>2$ classes while keeping as most discriminant information in the lower-dimensional space as possible.

Let $\boldsymbol{d} \in {\mathbb{R}}^{M}$ serve as the projection direction. In Fisher's discriminant analysis \citep{fisher1936use}, the Fisher criterion reads
\begin{equation}
{\underset{\boldsymbol{d}}{\max}} \frac{\boldsymbol{d}^{\top} \boldsymbol{S}_{\rm B} \boldsymbol{d}}{\boldsymbol{d}^{\top} \boldsymbol{S}_{\rm W} \boldsymbol{d}}, \label{fisher_criterion}
\end{equation}
which can be addressed by solving 
\begin{equation}
\begin{array}{l}
\boldsymbol{S}_{\rm B} \boldsymbol{d}=\lambda \boldsymbol{S}_{\rm W} \boldsymbol{d} \label{Eqn:generalized},
\end{array}
\end{equation}
where $\lambda$ is the Lagrange multiplier \citep{boyd2004convex}. This is also known as the generalized eigenvalue problem regarding $\boldsymbol{S}_{\rm B}$ and $\boldsymbol{S}_{\rm W}$, and $\boldsymbol{d}$ is the eigenvector corresponding to the non-zero eigenvalue ($\lambda$) in this situation. Then the transformation matrix can be formed by stacking up the $(C-1)$ eigenvectors corresponding to the $(C-1)$ largest eigenvalues in Eq. \eqref{Eqn:generalized}.
When the number of samples $N$ is small and/or the dimensionality of the data $M$ is big, $\boldsymbol{S}_{\rm W}$ is generally singular in practice. 
This could be dealt with by adding a small perturbation on $\boldsymbol{S}_{\rm W}$, e.g.,
\begin{equation}
\hat{\boldsymbol{S}}_{\rm W}=\boldsymbol{S}_{\rm W}+\delta \boldsymbol{I},
\end{equation}
where $\boldsymbol{I}$ is the identity matrix and $\delta$ is a relatively small value (e.g., $5\times 10^{-3}$) such that $\hat{\boldsymbol{S}}_{\rm W}$ is therefore invertible. 
The discriminant directions can then obtained by conducting the eigenvalue decomposition of $\hat{\boldsymbol{S}}_{\rm W}^{-1} \boldsymbol{S}_{\rm B}$ and finding the $(C-1)$ eigenvectors corresponding to the $(C-1)$ largest eigenvalues.

\subsubsection{{Binary Classification Problem}} \label{Sec:da-2c}
Different from Fisher criterion given in Eq. \eqref{fisher_criterion}, which can only produce one discriminant direction in the binary classification scenario, the method proposed in \citep{foley1975optimal} can discover more discriminant directions. It is optimal in the sense that a set of projection directions $\{\boldsymbol{d}_{k}\}_{k=1}^n$ is determined under a variety of constraints, see details below. 

The Fisher criterion ({\em cf.} Eq. \eqref{fisher_criterion}) for the binary classification problem reads
\begin{equation}
{\cal R}(\boldsymbol{d}) = \frac{\boldsymbol{d}^{\top} \tilde{\boldsymbol{S}}_{\rm B} \boldsymbol{d}}{\boldsymbol{d}^{\top} \tilde{\boldsymbol{S}}_{\rm W} \boldsymbol{d}}.
\label{Eqn:da-fc}
\end{equation}
Note that ${\cal R}(\boldsymbol{d})$ is independent of the magnitude of $\boldsymbol{d}$. The first discriminant direction $\boldsymbol{d}_1$ is discovered by maximising ${\cal R}(\boldsymbol{d})$, and then we have 
\begin{equation}
\boldsymbol{d}_{1} = \alpha_1 {\tilde{\boldsymbol{S}}_{\rm W}}^{-1}\boldsymbol{s}_{\rm b},
\label{eqn:da-d1}
\end{equation}
where $\alpha_1$ (i.e., $\alpha_{1}^{2}=(\boldsymbol{s}_{\rm b}^{\top}[\tilde{\boldsymbol{S}}_{\rm W}^{-1}]^{2} \boldsymbol{s}_{\rm b})^{-1}$) is the normalising constant such that $\|\boldsymbol{d}_{1}\|_2 = 1$ (and recall $\boldsymbol{s}_{\rm b}$ is the difference of the means of the two classes).
The second discriminant direction $\boldsymbol{d}_{2}$ is required to maximise ${\cal R}(\boldsymbol{d})$ in Eq. \eqref{Eqn:da-fc} and be orthogonal to $\boldsymbol{d}_{1}$. It can be found by the method of Lagrange multipliers, i.e., finding the stationary points of 
\begin{equation}
{\cal R}(\boldsymbol{d}_2) - \lambda[\boldsymbol{d}_{2}^{\top} \boldsymbol{d}_{1}],
\label{Eqn:d2}
\end{equation}
where $\lambda$ is the Lagrange multiplier. We can then obtain
\begin{equation}
\boldsymbol{d}_{2} = 
\alpha_2 \left(\tilde{\boldsymbol{S}}_{\rm W}^{-1} - \frac{\boldsymbol{s}_{\rm b}^{\top}({\tilde{\boldsymbol{S}}_{\rm W}}^{-1})^{2} \boldsymbol{\boldsymbol{s}_{\rm b}}}{\boldsymbol{s}_{\rm b}^{\top}({\tilde{\boldsymbol{S}}_{\rm W}}^{-1})^{3} \boldsymbol{s}_{\rm b}} ({\tilde{\boldsymbol{S}}_{\rm W}}^{-1})^{2}\right) {\boldsymbol{s}_{\rm b}},
\end{equation}
where $\alpha_2$ is the normalising constant such that $\|\boldsymbol{d}_{2}\|_2 = 1$.

The above procedure can be extended to any number of directions (until the number of features $M$) recursively as follows. The $n$-th discriminant direction $\boldsymbol{d}_{n}$ is required to maximise ${\cal R}(\boldsymbol{d})$ in Eq. \eqref{Eqn:da-fc} and be orthogonal to $\boldsymbol{d}_{k}, k = 1,2, \ldots, n-1$.
It can be shown that
\begin{equation}
\label{dn}
\boldsymbol{d}_{n} = \alpha_n
{\tilde{\boldsymbol{S}}_{\rm W}}^{-1}\left\{\boldsymbol{s}_{\rm b}-\left[\boldsymbol{d}_{1} \cdots \boldsymbol{d}_{n-1}\right] \boldsymbol{S}_{n-1}^{-1}\left[\begin{array}{c}
1 / \alpha_{1} \\
0 \\
\vdots \\
\vdots \\
0
\end{array}\right]\right\},
\end{equation}
where $\alpha_n$ is the normalising constant such that $\|\boldsymbol{d}_{n}\|_2 = 1$ and the $(i, j)$ entries of $\boldsymbol{S}_{n-1} \in \mathbb{R}^{(n-1)\times (n-1)}$ are defined as 
\begin{equation} \label{eqn-da-sn}
\boldsymbol{d}_{i}^\top {\tilde{\boldsymbol{S}}_{\rm W}}^{-1} \boldsymbol{d}_{j}, \ \  1\le i, j \le n-1.
\end{equation}
The whole procedure of finding $L$ number of discriminant vectors $\{\boldsymbol{d}_{n}\}_{n=1}^L$ is summarised in Algorithm \ref{alg:da-fs}.

\begin{algorithm}[t]
\caption{LDA for binary classification}\label{alg:da-fs}
\begin{algorithmic}
\STATE \textbf{Require:}
$\{\boldsymbol{y}_i\}_{i=1}^N$ and $L\le M$, i.e., the given samples and the number of discriminant vectors.
\STATE \hspace{0.3cm} \textbf{Compute}
$\tilde{\boldsymbol{S}}_{\rm W}$ and $\boldsymbol{s}_{\rm b}$ in Eq. \eqref{eqn:bi-intra-inter};
\STATE \hspace{0.3cm} \textbf{Compute} $\boldsymbol{d}_{1}$ using Eq. \eqref{eqn:da-d1} and $\boldsymbol{S}_{1}$ using Eq. \eqref{eqn-da-sn};
\STATE \hspace{0.3cm} $n$ = 1;
\STATE \hspace{0.3cm} \textbf{for} {$n < L$} \textbf{do}
\STATE \hspace{0.8cm} ${n}={n}+1$;
\STATE \hspace{0.8cm} \textbf{Compute} $\boldsymbol{d}_{n}$ using Eq. \eqref{dn};
\STATE \hspace{0.8cm} \textbf{Compute} $\boldsymbol{S}_{n}$ using Eq. \eqref{eqn-da-sn};
\STATE \hspace{0.3cm} \textbf{end for}
\STATE  \textbf{Return} $\{\boldsymbol{d}_{n}\}_{n=1}^L$
\end{algorithmic}
\label{alg1}
\end{algorithm}

Similar to how each singular vector correlates to a singular value, each discriminant vector $\boldsymbol{d}_{n}$ corresponds to a “discrim-value” say $\gamma_{n}$, where
\begin{equation}
\gamma_{n}=\frac{\boldsymbol{d}_n^{\top} \tilde{\boldsymbol{S}}_{\rm B} \boldsymbol{d}_n}{\boldsymbol{d}_{n}^{\top} \tilde{\boldsymbol{S}}_{\rm W} \boldsymbol{d}_{n}}.
\label{Eqn:dv}
\end{equation}
The discriminant vectors $\{\boldsymbol{d}_{n}\}_{n=1}^L$ are naturally ordered by their discriminant values, following {$\gamma_{1} \geq \gamma_{2} \geq \cdots \geq \gamma_{L} \geq 0$}. 

The DA subspace formed by $\{\boldsymbol{d}_{n}\}_{n=1}^L$ offers considerable potential for feature extraction and dimensionality reduction in many fields like pattern recognition \citep{foley1975optimal}.
For example, face recognition has been enhanced by LDA \citep{chelali2009linear} outperforming PCA in many cases.

\subsection{Non-negative Matrix Factorization}
In the process of matrix factorization, reconstructing a low-rank approximation for the data matrix $\boldsymbol{Y}$ is of great importance.
NMF is a technique dealing with $\boldsymbol{Y} \geq 0$ whose entries are all non-negative \citep{NIPS2000_f9d11525}, with great achievements in many fields such as signal processing \citep{dong2021transferred}, biomedical engineering \citep{leuschner2019supervised}, pattern recognition \citep{lee1999learning} and image processing \citep{chen2021deep}. The sparsity of the NMF subspace has also received extensive attention. In genomics, for example, the work in \citep{Brunet4164} factorized gene expression matrices across different experimental conditions, showing that the sparsity of NMF contributes to decreasing noise and extracting biologically meaningful features.
The purpose of NMF is to find two non-negative and low-rank matrices, i.e., one base matrix $\boldsymbol{X} \in \mathbb{R}^{p\times M}$ and one coefficient matrix $\boldsymbol{K} \in \mathbb{R}^{N\times p}$, satisfying
\begin{equation}
\boldsymbol{Y} \approx \boldsymbol{K} \boldsymbol{X},
\label{eqn:nmf}
\end{equation}
where $p < \min \{M, N\}$. Let $\boldsymbol{K} = (\boldsymbol{k}_1, \boldsymbol{k}_2, \cdots, \boldsymbol{k}_N)^\top$. We have $\boldsymbol{y}_i^\top \approx \boldsymbol{k}_i^\top \boldsymbol{X}, 1\le i \le N$. In other words, every sample $\boldsymbol{y}_i$ can be represented by a linear combination of the rows of $\boldsymbol{X}$ with the components in $\boldsymbol{k}_i$ serving as weights. Therefore, $\boldsymbol{X}$ is also known as consisting of basis vectors which can project the data matrix $\boldsymbol{Y}$ into a low-dimensional subspace. The number of basis vectors $p$ will affect the degree of approximation to the data matrix $\boldsymbol{Y}$.

Finding $\boldsymbol{K}$ and $\boldsymbol{X}$ satisfying Eq. \eqref{eqn:nmf} can be addressed by solving the following minimisation problem:
\begin{equation} \label{eqn:nmf-min}
\min_{\boldsymbol{K, X}}\|\boldsymbol{Y-K X}\|_{F}^{2}, \quad
\text {s.t.} \ \boldsymbol{K} \geq 0, \boldsymbol{X} \geq 0,
\end{equation}
where $\|\cdot\|_{F^{}}$ is the Frobenius norm.
To solve problem \eqref{eqn:nmf-min}, a common technique is to update $\boldsymbol{K}$ and $\boldsymbol{X}$ alternatively, i.e.,
\begin{equation}
\label{eqn:nmf-update}
\boldsymbol{K} \leftarrow \boldsymbol{K} \odot \frac{\boldsymbol{Y X}^{\top}}{\boldsymbol{K X X}^{\top}}, \quad
\boldsymbol{X} \leftarrow \boldsymbol{X} \odot \frac{\boldsymbol{K}^{\top} \boldsymbol{Y}}{\boldsymbol{K}^{\top} \boldsymbol{K X}},
\end{equation}
where $\odot$ denotes the pointwise product. For more algorithmic details please refer to e.g. \citep{NIPS2000_f9d11525}.

NMF is an unsupervised method that decomposes the data matrix without utilising the class label information.
Regarding the binary classification problem, the SNMF (supervised NMF) proposed in \citep{leuschner2019supervised} extends the standard unsupervised NMF approach by exploiting feature extraction and integrating the cost function of the classification method into NMF. In SNMF, the classification labels are incorporated into the algorithms to extract the specific data patterns relevant to the respective classes. 
The whole algorithm of SNMF is provided in Appendix \ref{Appendix-snmf}.

\begin{figure*}[!t]
    \centering
    \includegraphics[width=6.6in]{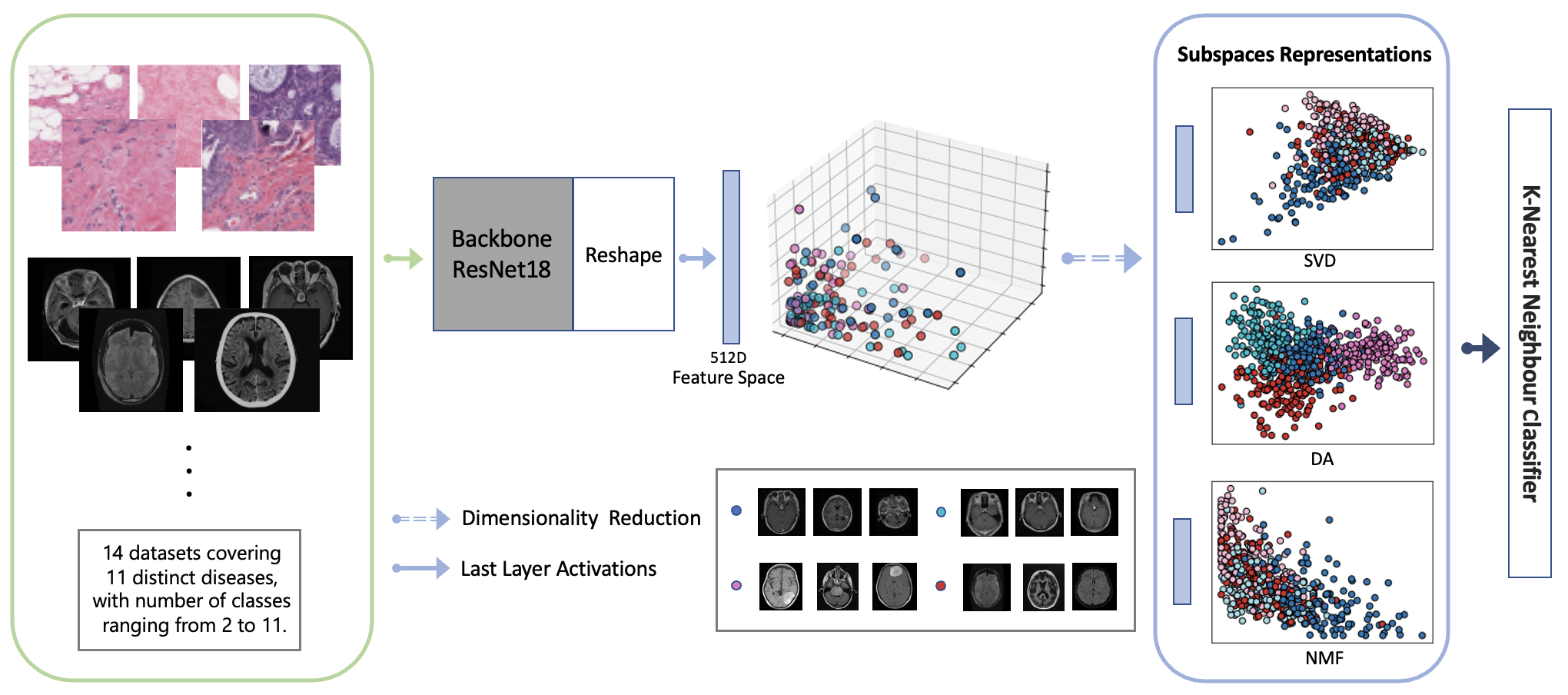}
    \centering
    \caption{Few-shot learning schematic diagram on different subspaces. From left to right: A pre-trained deep neural network (e.g. ResNet18) to solve a large natural image classification problem is exploited to extract features of medical images (i.e., inputs in the green box), and then the extracted features are projected to subspace representations (i.e., outputs in the blue box), followed by a classifier (e.g. KNN) delivering the classification results. The extracted features for individual images are visualised as dots with different colours representing different classes (i.e., the middle of the diagram).}
    \label{fig_sim}
\end{figure*}

\section{Few-shot Learning on Subspace Representations}
\label{Sec:fl-art}
We deploy few-shot learning techniques in investigating medical imaging particularly in the data scarcity scenario. We consider problems in which the feature space dimensionality is usually high in comparison to the number of images we have; hence subspace representations are sought. The adopted few-shot learning scheme on subspace feature representations and experimental settings are presented in what follows.

\subsection{Framework }\label{subsec:expset-fl}
The deployed and enhanced few-shot learning schematic diagram on different subspaces is shown in Figure \ref{fig_sim}. Firstly, a pre-trained deep neural network (e.g. ResNet18) to solve a large natural image classification problem is prepared and then used to extract features of medical images in the given datasets (i.e., the green box in Figure \ref{fig_sim}). After that the extracted features are projected to subspace representations (i.e., the blue box in Figure \ref{fig_sim}).  In this paper, we consider three different methods (i.e., SVD, DA and NMF) described in Section \ref{Sec:subspace-rep} to achieve this. Finally, a classifier (e.g. the $K$-Nearest Neighbour (KNN) or Support Vector Machine (SVM)) is employed  to perform few-shot learning -- predicting the final classification results. 
Extensive exploration in terms of the benefits of different subspace representations and insightful suggestions and comparisons in the regime of few-shot learning in medical imaging will be conducted in Section \ref{sec:results}. 

\subsection{Experimental Settings}\label{subsec:exp-set}
We explore 14 datasets covering 11 distinct diseases, with the number of classes ranging from 2 to 11, see Section \ref{sec:data} for more detail. The pre-trained deep model, ResNet18, is used as the source model in our experiment. Each input is pre-processed and pixel-wised by subtracting its mean and being divided by the standard deviation without data augmentation. The feature space is from the features in the penultimate layer of the pre-trained model (ResNet-18) extracted by PyTorch hooks \citep{pytorch}, yielding a $512$-dimensional feature vector for each image. The low-dimensional representations are then generated from the introduced methods. The number of iterations related to NMF and SNMF is set to $3000$ to ensure convergence. The mean result of the KNN classifier with selected $K$ (with values of 1, 5, 10 and 15) nearest neighbours is used to evaluate the final performance. Except for KNN, we also implement SVM as the classifier for comparison. The detailed experimental setting and results of the SVM classifier are shown in Appendix \ref{svm-nmf-da}. To quantify the uncertainty of the classification accuracy and produce more reliable quantitative results, we present averages and standard deviations across 10 distinct times of random samplings in each dataset. In addition to the accuracy, the reconstruction error of NMF at different random initialization is conducted to demonstrate its convergence. Moreover, we also compare our method with other well-known few-shot learning algorithms like the prototypical network \citep{snell2017prototypical}. The experimental setup and results are presented in Appendix \ref{Prototypical}.

\section{Data}\label{sec:data}
A total of $14$ different datasets covering a range of problems in diagnostics are employed for our
empirical work. The number of classes ranges from 2 to 11 and the imaging modalities include X-rays, CT scans, MRI and Microscope. The datasets with {\tt MNIST}
within their name come from a benchmark family referred to as {\tt MedMNIST}\footnote{\url{https://medmnist.com/}}.
In order to illustrate the regime of few-shot learning, randomly sampled subsets of the whole individual datasets are used for our training and test. The corresponding data split for each class in the training and test sets for all the datasets is presented in Table \ref{table:data_description}.
It is worth noting that our intention is not to compare with previously published results which have used the whole individual datasets.
For ease of reference, brief descriptions of these individual datasets together with our implementations are given below. 

1) {\tt BreastCancer (IDC)} data \citep{breastcancer,cruz2014automatic}  is a binary classification problem sampled from digitised whole slide histopathology images. The source of the data is from $162$ women diagnosed with Invasive Ductal Carcinoma (IDC), the most common phenotypic subtype in breast cancers. From these annotated images $277,524$ patches had been segmented. An accuracy of $84.23\%$ using the whole dataset is reported in \citep{cruz2014automatic}. 

2) {\tt BrainTumor} data \citep{cheng_2017,cheng2015enhanced} is a four-category problem,  consisting of $7,022$ images of human brain MRI images, three types of tumours (i.e., glioma, meningioma and pituitary), and a control group. 

3) {\tt CovidCT} data \citep{he2020sample} is a binary classification problem, which is of great interest due to the COVID-19 pandemic. It contains $349$ CT scans that are positive for COVID-19 and $397$ negatives that are normal or contain other types of diseases. Two-dimensional slices from the scans are used in the study.

4) {\tt DeepDRiD} data \citep{DeepDRi} is a five-category problem. Diabetic retinopathy is a prevalent eyesight condition in eye care. With early detection and treatment, the majority of these disorders may be controlled or cured. In this dataset, a total of $2,000$ regular fundus images were acquired using Topcon digital fundus camera from 500 patients. 

5) {\tt BloodMNIST} data \citep{acevedo2020dataset} is an eight-category problem, including a total of $17,092$ images. It consists of individual normal cells, captured from individuals without infection, hematologic or oncologic disease and free of any pharmacologic treatment at the time of blood collection. 

6) {\tt BreastMNIST} data \citep{al2020dataset} is a binary classification problem, including a total of $780$ breast ultrasound images. An accuracy of $94.62\%$ is claimed in \citep{moon2020computer} in the computer-aided diagnostic (CAD) setting on the whole dataset. The grayscale images are replicated in order to match the pre-trained model. 

7) {\tt DermaMNIST} data \citep{tschandl2018ham10000,codella2019skin} is a multi-source dermatoscopic image collection of common pigmented skin lesions. It contains $10,015$ dermatoscopic images, which are classified into seven diseases. 

8) {\tt OCTMNIST} data \citep{kermany2018identifying} is for retinal diseases, including a total of $109,309$ valid optical coherence tomography images, with four diagnostic categories. 

9) {\tt OrganAMNIST}, {\tt OrganCMNIST} and {\tt OrganSMNIST} datasets \citep{bilic2019liver} are eleven-category problem. They are benchmarks for segmenting liver tumours from 3D computed tomography images (LiTS). Organ labels were obtained using boundary box annotations of the 11 bodily organs studied, which are renamed from Axial, Coronal and Sagittal for  simplicity. Grayscale images were converted into RGB images through the instruction in \citep{medmnistv2}. 

10) {\tt PathMNIST} data \citep{kather2019predicting} is based on the study of using colorectal cancer histology slides to predict survival, including a total of 107,180 images and nine different types of tissues. An accuracy of $94\%$ was achieved in \citep{kather2019predicting} by training a CNN using transfer learning on a set of $7,180$ images from $25$ CRC patients. 

11) {\tt PneumoniaMNIST} data \citep{kermany2018identifying} is to classify pneumonia into two categories -- severe and mild. It consists of
$5,856$ paediatric chest X-ray images. The source images are grayscale, which are converted to RGB for training in the same manner as the {\tt OrganAMNIST} dataset.

12) {\tt TissueMNIST} data \citep{woloshuk2021situ} is derived from the Broad Bioimage Benchmark Collection. It consists of $236,386$ human kidney cortex cells, segmented and labelled into eight categories. An accuracy of $80.26\%$ was achieved in \citep{woloshuk2021situ} using a custom 3D CNN on the whole dataset.

\begin{table}[b]
\caption{Data split for each class in the training and test sets of each dataset.}
\label{table:data_description}
\centering
\setlength{\tabcolsep}{9mm}
\renewcommand\arraystretch{1.1}
\begin{tabular}{c||c|c|c}
\toprule 
\multirow{2}{*}{Datasets}  & \multicolumn{1}{c|}{\multirow{2}{*}{\#Classes}} & \multicolumn{2}{c}{\#Samples for each class} \\ 
\cline{3-4} & \multicolumn{1}{c|}{} & Training & Test     \\ 
\midrule 
{\tt BreastCancer}\citep{breastcancer}             & 2       & 300      & 40   \\
{\tt BrainTumor}\citep{cheng_2017}                & 4       & 160      & 40   \\
{\tt CovidCT}\citep{he2020sample}                  & 2       & 300      & 40   \\
{\tt DeepDRiD} \citep{DeepDRi}                     & 5       & 118      & 29   \\
{\tt BloodMNIST}\citep{acevedo2020dataset}         & 8       & 75       & 25   \\
{\tt BreastMNIST}\citep{al2020dataset}             & 2       & 263      & 88   \\
{\tt DermaMNIST}\citep{tschandl2018ham10000}       & 7       & 75       & 25   \\
{\tt OCTMNIST}\citep{kermany2018identifying}       & 4       & 150      & 50   \\
{\tt OrganAMNIST}\citep{bilic2019liver}            & 11      & 50       & 15   \\
{\tt OrganCMNIST}\citep{bilic2019liver}            & 11      & 50       & 15   \\
{\tt OrganSMNIST}\citep{bilic2019liver}            & 11      & 50       & 15   \\
{\tt PathMNIST}\citep{kather2019predicting}        & 9       & 60       & 20   \\
{\tt PneumoniaMNIST}\citep{kermany2018identifying} & 2       & 262      & 87   \\
{\tt TissueMNIST}\citep{woloshuk2021situ}          & 8       & 65       & 20   \\ 
\midrule 
\end{tabular}
\end{table}

\section{Experimental Results}\label{sec:results}
In this section, we investigate the performance of the few-shot learning scheme described in Section \ref{Sec:fl-art} on subspace representations using SVD, DA and NMF. Note, importantly, that our main interest is to introduce DA and NMF as alternative subspace representations to SVD in the regime of few-shot learning in medical imaging. 
In addition to the comparison between the SVD, DA and NMF subspaces, we  also compare them with other relevant feature selection, dimensionality reduction, and few-shot learning methods.
For visual inspection, we visualise the subspace distributions of SVD, DA and NMF by T-SNE built-in function in Python (see the results in Appendix \ref{TSNE visualisation}).

\begin{table}[!t]
\caption{Few-shot learning classification accuracy on 14 medical datasets with the KNN classifier. The original feature space and the subspaces obtained by SVD and DA are tested. The number of classes ranging from 2--11 shows the multi-modal structure in each dataset.}
\label{table:acc}
\centering
\setlength{\tabcolsep}{1mm}
\begin{tabular}{c||c|c|c|c}
\toprule
Datasets & Feature Space & SVD & DA & Classes
\\
\midrule \midrule
{\tt BreastCancer}\citep{breastcancer} &63.25$\pm$4.80&60.85$\pm$8.44&\textbf{66.76$\pm$5.39}&2
\\
{\tt BrainTumor}\citep{cheng_2017}
&68.63$\pm$4.17&54.63$\pm$4.73&\textbf{73.19$\pm$2.98}&4
\\
{\tt CovidCT}\citep{he2020sample}
&\textbf{77.11$\pm$2.89}&66.68$\pm$6.82&70.58$\pm$3.39&2
\\
{\tt DeepDRiD} \citep{DeepDRi}
&48.25$\pm$6.94&36.19$\pm$4.98&\textbf{55.11$\pm$4.19}&5
\\
{\tt BloodMNIST}\citep{acevedo2020dataset}
&37.49$\pm$3.88&37.10$\pm$3.91&\textbf{54.33$\pm$3.56}&8
\\
{\tt BreastMNIST}\citep{al2020dataset}
&69.78$\pm$3.79&68.45$\pm$4.19&\textbf{70.08$\pm$3.58}&2
\\
{\tt DermaMNIST}\citep{tschandl2018ham10000}
&25.03$\pm$4.64&21.16$\pm$3.29&\textbf{33.52$\pm$3.13}&7
\\
{\tt OCTMNIST}\citep{kermany2018identifying}
&31.61$\pm$4.44&28.25$\pm$3.60&\textbf{34.85$\pm$3.13}&4
\\
{\tt OrganAMNIST}\citep{bilic2019liver}
&32.65$\pm$2.58&35.30$\pm$3.26&\textbf{49.67$\pm$2.98}&11
\\
{\tt OrganCMNIST}\citep{bilic2019liver}
&25.80$\pm$3.09&26.88$\pm$3.13&\textbf{45.93$\pm$4.14}&11
\\
{\tt OrganSMNIST}\citep{bilic2019liver}
&24.80$\pm$2.37&24.18$\pm$2.54&\textbf{39.59$\pm$2.64}&11
\\
{\tt PathMNIST}\citep{kather2019predicting}
&33.97$\pm$2.37&38.47$\pm$4.59&\textbf{58.68$\pm$3.75}&9
\\
{\tt PneumoniaMNIST}\citep{kermany2018identifying}
&70.43$\pm$3.70&61.60$\pm$7.48&\textbf{73.76$\pm$5.22}&2
\\
{\tt TissueMNIST}\citep{woloshuk2021situ}
&18.89$\pm$2.80&16.88$\pm$2.42&\textbf{21.86$\pm$2.15}&8\\
\midrule
\end{tabular}
\end{table}

\subsection{Discriminant versus Principal Component Subspaces}
We first conduct comparison between DA and PCA.
Table \ref{table:acc} shows the few-shot learning classification accuracy on the 14 datasets/problems, comparing the feature space in its original dimension of the ResNet18 with the PCA and DA subspaces. The accuracy results are the average of $K$ values of KNN classifier chosen to be 1, 5, 10 and 15. We note that with a single exception of the {\tt CovidCT} dataset, principal component dimensionality reduction loses information about class separation, whereas the discriminant subspace representation maintains the separation extremely well, thereby showing significant improvement over the original feature space. In detail, in $11$ of the $14$ problems, the SVD subspace performs worse than the original feature space. In contrast, the DA subspace shows significant improvement over the corresponding SVD subspace in all the 14 problems; and in $13$ of the $14$ problems, the DA subspace shows significant improvement over the original feature space. Furthermore, Z-test was also carried out and it is confirmed that the results are statistically significant at $P$ values smaller than $10^{-3}$. 

\begin{figure}[tp]
  \centering
  \begin{center}
  \begin{minipage}[c]{0.4\textwidth}
    \centering
    \subfigure[{\tt PneumoniaMNIST} dataset]{\includegraphics[width=\textwidth]{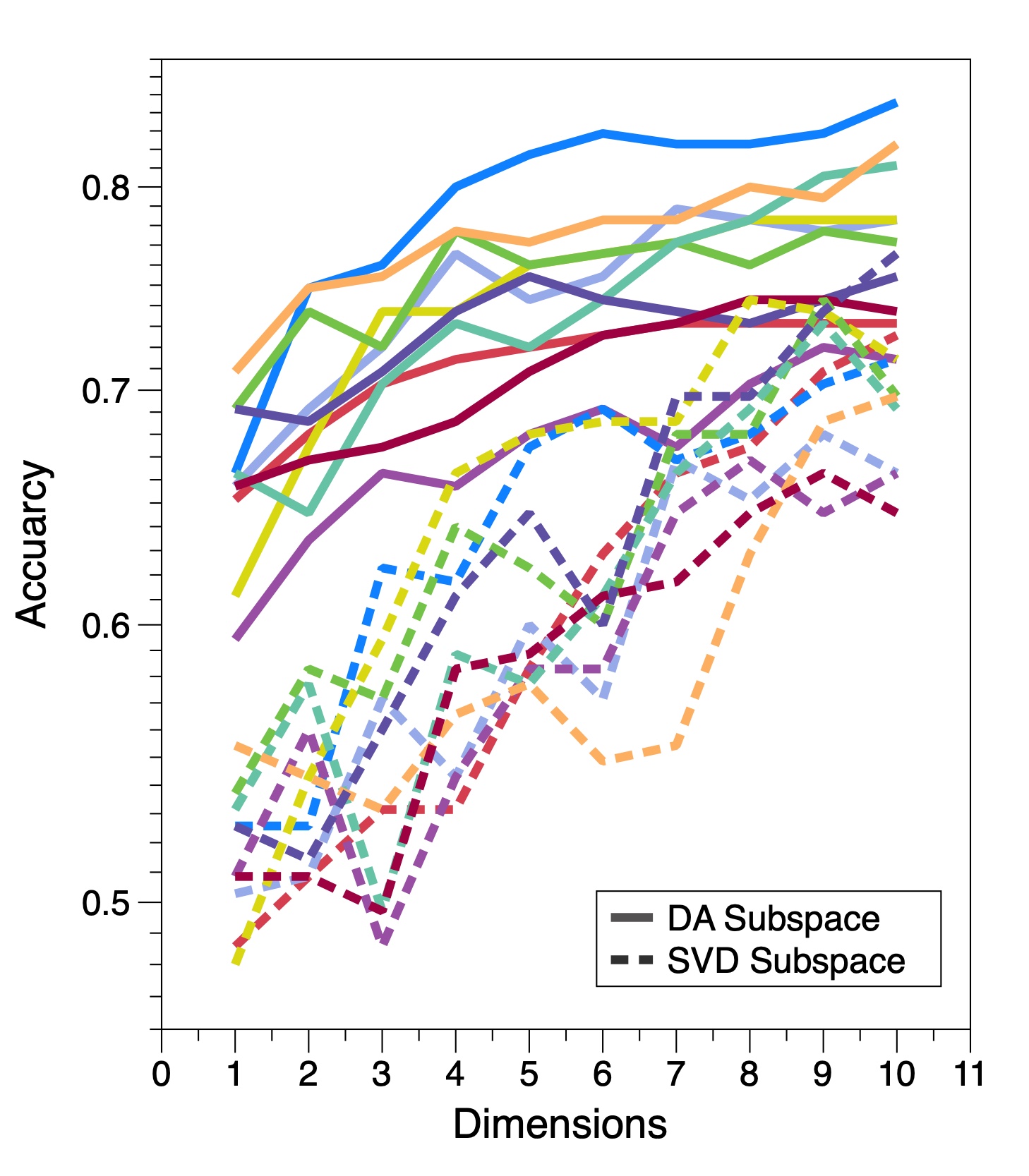}
  \label{fig:2classes}}
  \end{minipage}
  \hfill
  \begin{minipage}[c]{0.55\textwidth}
    \centering
    \subfigure[{\tt OrganAMNIST} dataset]{\includegraphics[width=\textwidth]{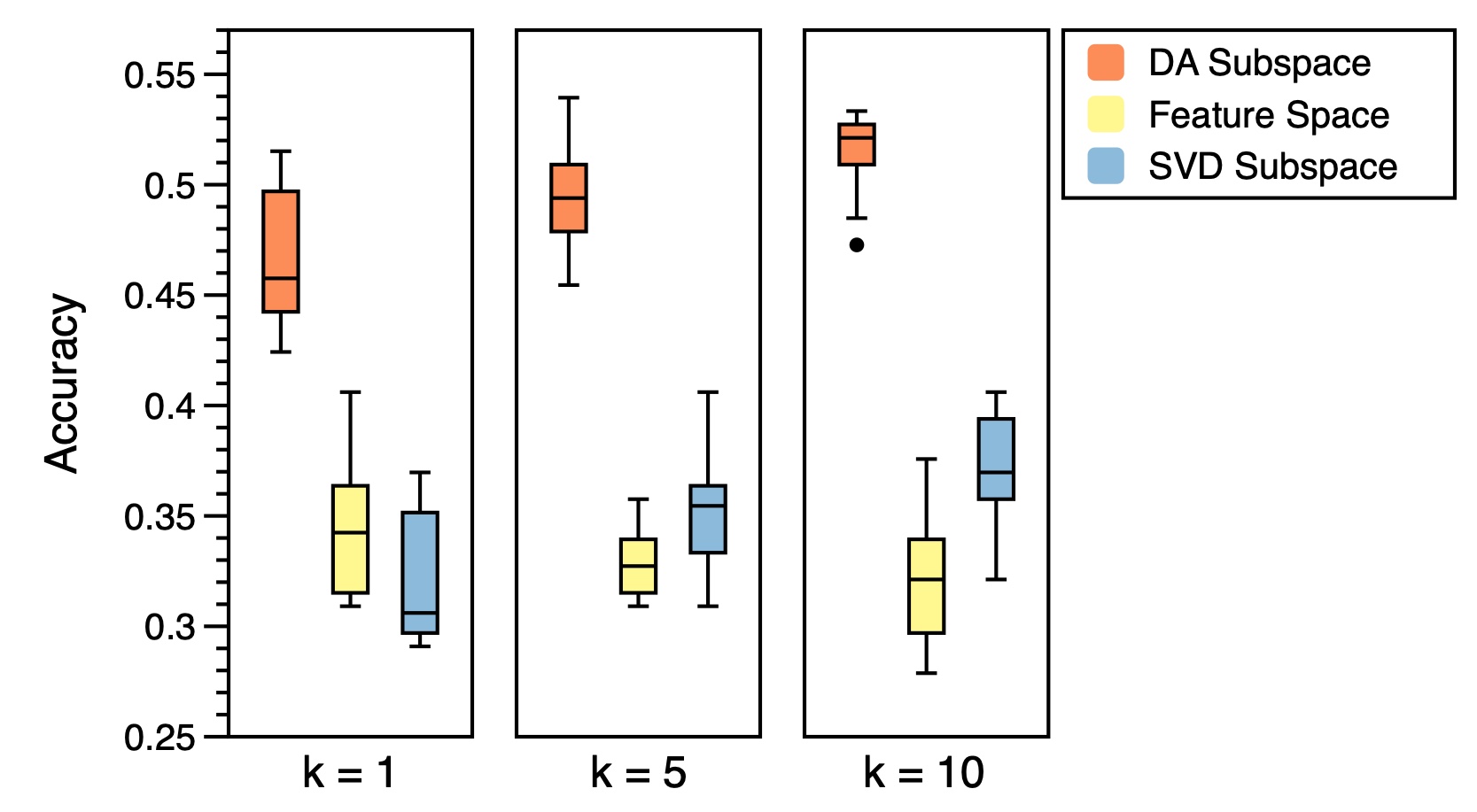}
  \label{fig:11classes}}
  \end{minipage}    
  \caption{Comparison between DA and PCA subspaces in terms of classification accuracy corresponding to different dimensions and different neighbourhood size $K$ in the KNN classifier. Figure (a) shows the DA subspace taken at different dimensions consistently outperform the SVD subspace ({\em cf.} Table \ref{table:acc} for the performance on the full $512$ dimensional feature space). 
  Figure (b) shows the excellent performance of the DA subspace against PCA and the original feature space, irrespective of the choice of $K$ in the classifier.}
  \end{center}
\end{figure}

We now evaluate the impact of the subspace dimensions on the classification accuracy for DA and SVD. Figure \ref{fig:2classes} shows how the classification accuracy varies as the subspace dimensions increase on the {\tt PneumoniaMNIST} dataset (consistent results are observed for other datasets). In particular, ten different random partitions of the training-test set are utilised to shuffle the data (which will make the results more credible) and dimensions from one to ten are investigated in Figure \ref{fig:2classes}.
We observe that the performance of both the DA and SVD methods increases monotonically corresponding to the number of dimensions, with the DA subspace consistently outperforming SVD. 
Given the performance achieved using the full set of features is $70.43 \pm 3.70$ in Table \ref{table:acc}, hence the increase for SVD is not sustainable beyond this point.

The effect of different neighbourhood size $K$ of the KNN classifier is reported in Figure \ref{fig:11classes}, where the eleven-class dataset {\tt OrganAMNIST} (consistent results are observed for other datasets) is used. Moreover, the performance of the SVD and DA subspaces with dimension equal to ten against the original feature space corresponding to $K = 1, 5$ and $10$ is evaluated in Figure \ref{fig:11classes}. Uncertainty in results is evaluated over $10$ random partitions of the training-test set, with $550$ and $165$ images for training and test, respectively. Figure \ref{fig:11classes} shows substantial improvement in DA subspace representation over both the original feature space and the SVD reduced subspace irrespective of the choice of $K$ in the KNN classifier.

\begin{figure}[!t]
\centering
\includegraphics[width=4in]{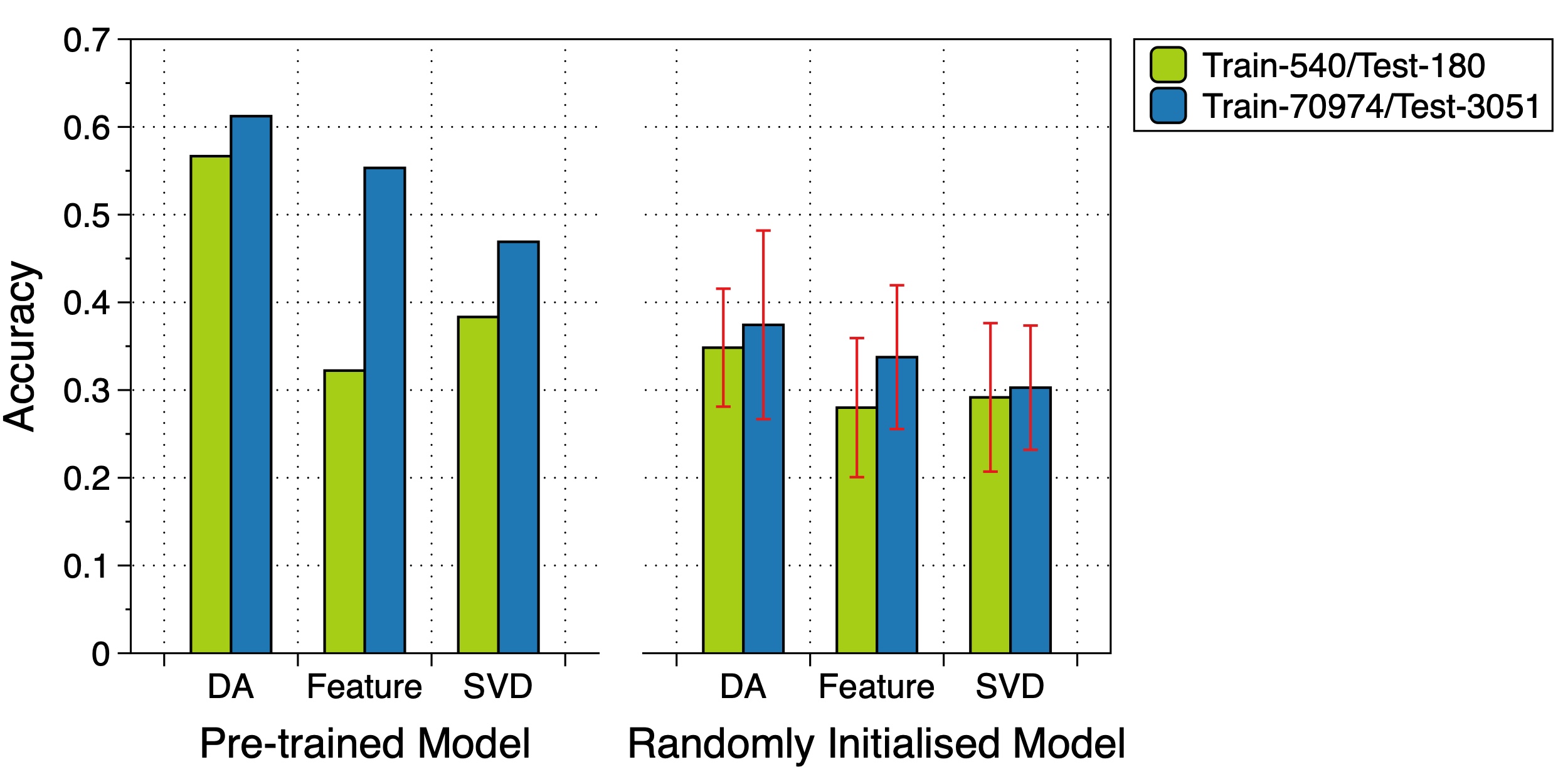}
\caption{Comparison between DA and PCA subspaces and the original feature space in terms of classification accuracy corresponding to different dataset sizes.
Dataset {\tt PathMNIST} with nine classes is used.
The left and right three pairs of bars in the panel are the results of the pre-trained model and the model with randomly initialised weights, respectively. The results reveal that the performance of the DA subspace always outperforms the SVD and the original feature space, irrespective of the choice of the data size. Moreover, the results achieved using the DA subspace are highly comparable to those obtained by using the entire dataset, whereas the results of SVD fall short.
}
\label{9classes}
\end{figure}

Finally, we investigate the effect of the dataset size on the performance of the methods compared.
Figure \ref{9classes} shows the results regarding the DA and PCA subspaces and the original feature space on a small subset (i.e., 540 and 180 images for training and test, respectively) of the dataset as well as the entire dataset (i.e., 70,974 and 3,051 images for training and test, respectively), where nine-class dataset {\tt PathMNIST} (consistent results are observed for other datasets) is used for illustration. The value $K$ in the KNN classifier is set to 5. In Figure \ref{9classes}, we also evaluate the effect of the pre-trained model on ImageNet versus the model whose weights are defined by random initialisation. The findings reveal that the performance of the DA subspace always outperforms the SVD and the original feature space, irrespective of the choice of the data size. Particularly, it also shows that, although utilising only 0.7\% of the entire dataset, the results achieved using the DA subspace are highly comparable to those obtained using the entire dataset, whereas the results of SVD fall short. This confirms that the DA subspace is more stable than the SVD subspace, providing a discriminative subspace ideal for classification problems. In passing, we also see that the performance of the pre-trained model is better than that of the model with randomly initialised weights, which fits our expectations. More results -- the comparison between DA and the manifold learning method Isomap (a non-linear dimensionality reduction process) -- on all the datasets are given in Appendix \ref{Appendix-nmf-result}.

\begin{table}[!t]
\centering
\caption{Binary classification accuracy comparison  between subspaces (i.e., SVD, NMF and SNMF) and the original feature space.}
\setlength{\tabcolsep}{1.3mm}{
\begin{tabular}{@{}c||cccc@{}}
\toprule
\multicolumn{1}{c||}{Datasets}& \multicolumn{1}{c|}{\begin{tabular}[c]{@{}c@{}}Feature\\  Space\end{tabular}} & \multicolumn{1}{c|}{SVD}                     & \multicolumn{1}{c|}{NMF}
& SNMF \\ \midrule \midrule
{\tt CovidCT}  & \multicolumn{1}{c|}{\textbf{77.11$\pm$2.89}} & \multicolumn{1}{c|}{75.75$\pm$2.82}          & \multicolumn{1}{c|}{74.96$\pm$1.98} & 76.47$\pm$2.25                                           \\
{\tt BreastCancer}  & \multicolumn{1}{c|}{63.25$\pm$4.80} & \multicolumn{1}{c|}{67.15$\pm$4.18}& \multicolumn{1}{c|}{68.75$\pm$4.45} & \textbf{69.53$\pm$4.98}                                   \\
{\tt PneumoniaMNIST} & \multicolumn{1}{c|}{70.43$\pm$3.70} & \multicolumn{1}{c|}{73.32$\pm$1.21} & \multicolumn{1}{c|}{75.08$\pm$2.07} & \textbf{75.28$\pm$2.02 }                                   \\
{\tt BreastMNIST}   & \multicolumn{1}{c|}{69.78$\pm$3.79} & \multicolumn{1}{c|}{70.92$\pm$2.88} & \multicolumn{1}{c|}{71.54$\pm$3.54} & \textbf{73.81$\pm$2.63}                                   \\ \bottomrule
\end{tabular}
}
\label{table:nmf-binary}
\end{table}

\begin{table}[!t]
\centering
\caption{Multiclass classification accuracy comparison  between subspaces (i.e., SVD and NMF) and the original feature space.}
\setlength{\tabcolsep}{4.0mm}{
\begin{tabular}{@{}c||c|c|c@{}}
\toprule
 Datasets &  Feature Space & SVD & NMF  \\ \midrule \midrule
{\tt DeepDRid} &  48.25$\pm$6.94  & 48.59$\pm$3.90 & \textbf{50.88$\pm$3.97} \\
{\tt BrainTumor} & 68.63$\pm$4.17   & 70.01$\pm$2.94 & \textbf{70.10$\pm$2.52} \\
{\tt BloodMNIST} & 37.49$\pm$3.88 & 44.88$\pm$1.59 & \textbf{45.18$\pm$2.30 }         \\
 {\tt DermaMNIST} & 25.03$\pm$4.64 & 28.96$\pm$2.44  & \textbf{29.79$\pm$2.01} \\
{\tt OCTMNIST} & 31.61$\pm$4.14  & 31.47$\pm$2.50 & \textbf{31.80$\pm$2.43}\\
{\tt OrganAMNIST}&  32.65$\pm$2.58 & 38.62$\pm$1.29  & \textbf{39.00$\pm$2.23} \\
{\tt OrganCMNIST}  &  25.80$\pm$3.09 & \textbf{32.92$\pm$2.77} & 32.11$\pm$2.07  \\
{\tt OrganSMNIST}&  24.80$\pm$2.37 & \textbf{29.12$\pm$1.46} & 29.01$\pm$1.51  \\
{\tt PathMINST}&  33.97$\pm$2.37 & 43.35$\pm$1.23 & \textbf{44.74$\pm$1.55} \\
{\tt TissueMNIST}&  18.89$\pm$2.80 & \textbf{21.03$\pm$1.58} & 20.75$\pm$1.52 \\
\bottomrule
\end{tabular}
}
\label{table:nmf-multi}
\end{table}

\subsection{Non-negative Matrix 
Factorization Subspace}
The classification accuracy of the NMF subspace (including NMF and SNMF) and the comparison with the SVD subspace and the original feature space on the binary class and multiclass problems are shown in Tables \ref{table:nmf-binary} and \ref{table:nmf-multi} respectively. The SNMF subspace is only limited to the binary class problem and the dimension of related subspaces is kept as $30$. 
It shows that, generally, the subspace representations (either SVD or NMF) deliver better performance than the original feature space. With SNMF marginally outperforming NMF in binary classification tests, NMF and SVD subspace both perform comparably and the trend is also preserved in multiclass classification problems. This prompts NMF can be a viable alternative to SVD, particularly when  sparse representation is of great interest.

Different from the dimension selected in the DA subspace, the dimension of the NMF/SVD subspaces is retained as 30. Mainly because the NMF (including SNMF) approximates the original data with the product of two matrices and is affected by the selected rank during decomposition. While for the DA subspace, the dimensions are determined by the number of classes for the multiclass problems.  Our results show that the performance of NMF is stable only after reaching a specific dimension, which is similar to the selection of the number of eigenvectors in SVD. Detailed trends regarding the performance of NMF and SVD subspaces on the 14 datasets against the changes in dimension are presented in Appendix \ref{Appendix-nmf-result}, including the comparison with the non-linear dimensionality reduction method Isomap in Figure \ref{all_results}.

\begin{figure}[htp]
  \centering
  \includegraphics[width=3.0in]{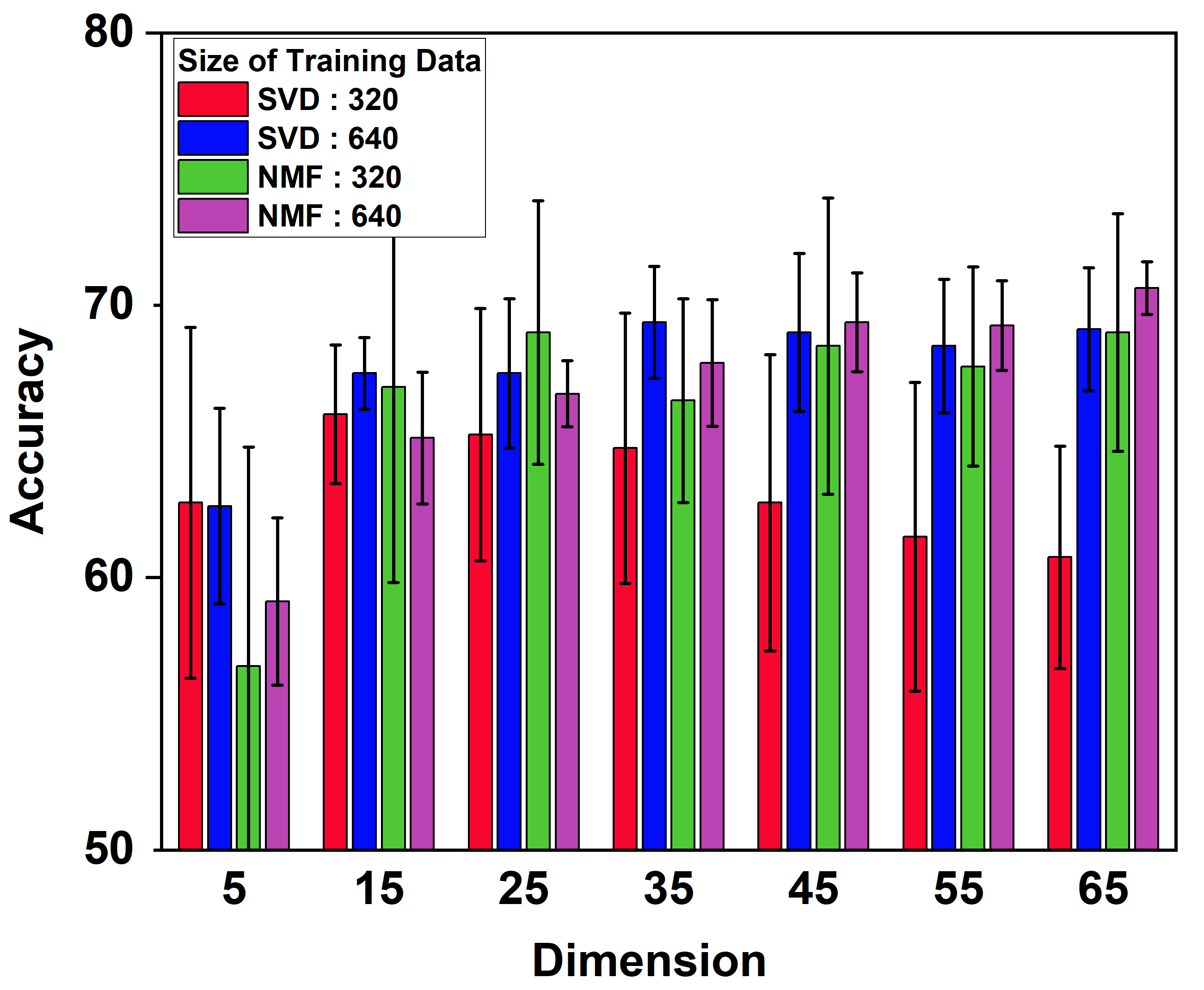}
  \caption{Comparison between NMF and SVD subspaces in terms of classification accuracy corresponding to different dataset size as the subspace dimension changes. Dataset {\tt BrainTumor} with four classes is used. Uncertainty is evaluated over 5 random partitions of the training-test set; and two types of training datasets with 320 and 640 images are created. The value $K=5$ is used in the KNN classifier. It shows that the performance of the NMF is stable for both types of datasets, whereas SVD suffers dimensional issues in the small dataset (with 320 images).}
  \label{datavolumn}
  \vspace{-0.15cm}
\end{figure}

\begin{figure}[!t]
  \centering
  \includegraphics[width=4.5in]{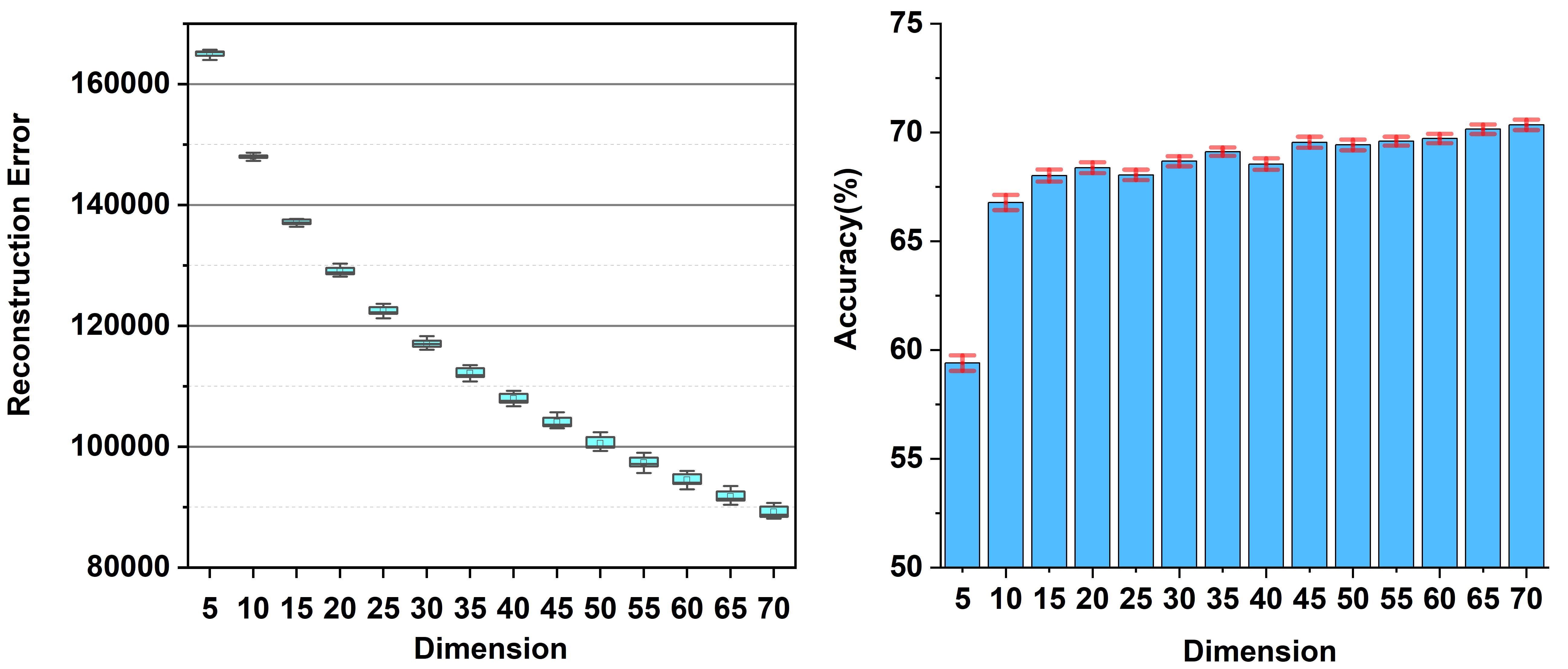}
\caption{Performance of the NMF subspace corresponding to different initialisation as the dimension of the subspace increases. The left and right panels respectively show the reconstruction error and the classification performance with $20$ random NMF initializations on the {\tt BrainTumor} dataset, indicating that the performance of NMF is quite stable corresponding to different dimensions with random initialisation except for the 5-dimensional subspace.
}
\label{reconstrution}
\vspace{-0.45cm}
\end{figure}

Additionally, we investigate the stability and uncertainty of NMF from the viewpoints of dataset size and the effects of random NMF initialization in various dimensions, respectively.
Figure \ref{datavolumn} describes how the volume of datasets influences the classification performance as the subspace dimensions ranging from 5 to 65 on the {\tt BrainTumor} dataset. Two training datasets with the size of 320 and 640 images are created for the SVD and NMF subspaces, represented by different colour bars. It shows that on the big dataset (with 640 images), SVD and NMF are quite similar (see the blue and purple bars).  On the small dataset (with 320 images), the NMF subspace outperforms the SVD subspace (see the red and green bars).  SVD suffers from dimension issues in the small dataset since it performs gradually worse rather than better when the dimension becomes higher (e.g. when the dimension increases from 15 to 65). In contrast, the results of the NMF subspace are relatively stable in different dimensions and have similar accuracy. Although NMF behaves not good in extremely low dimensions (such as 5 dimensions), it gets improved as the dimension increases, which is consistent with the statement mentioned before.  The uncertainty of NMF is evaluated by randomly initialising the NMF corresponding to different dimensions. In Figure \ref{reconstrution}, the left and right images show the reconstruction error and the classification performance with 20 random NMF initializations on the {\tt BrainTumor} dataset. It reveals that the reconstruction error decreases as the dimensionality increases and the performance of NMF is quite stable corresponding to different dimensions with random initialisation.

\subsection{Role of the Feature Extractor}

\begin{figure}[htp]
\centering
\includegraphics[width=3.5in]{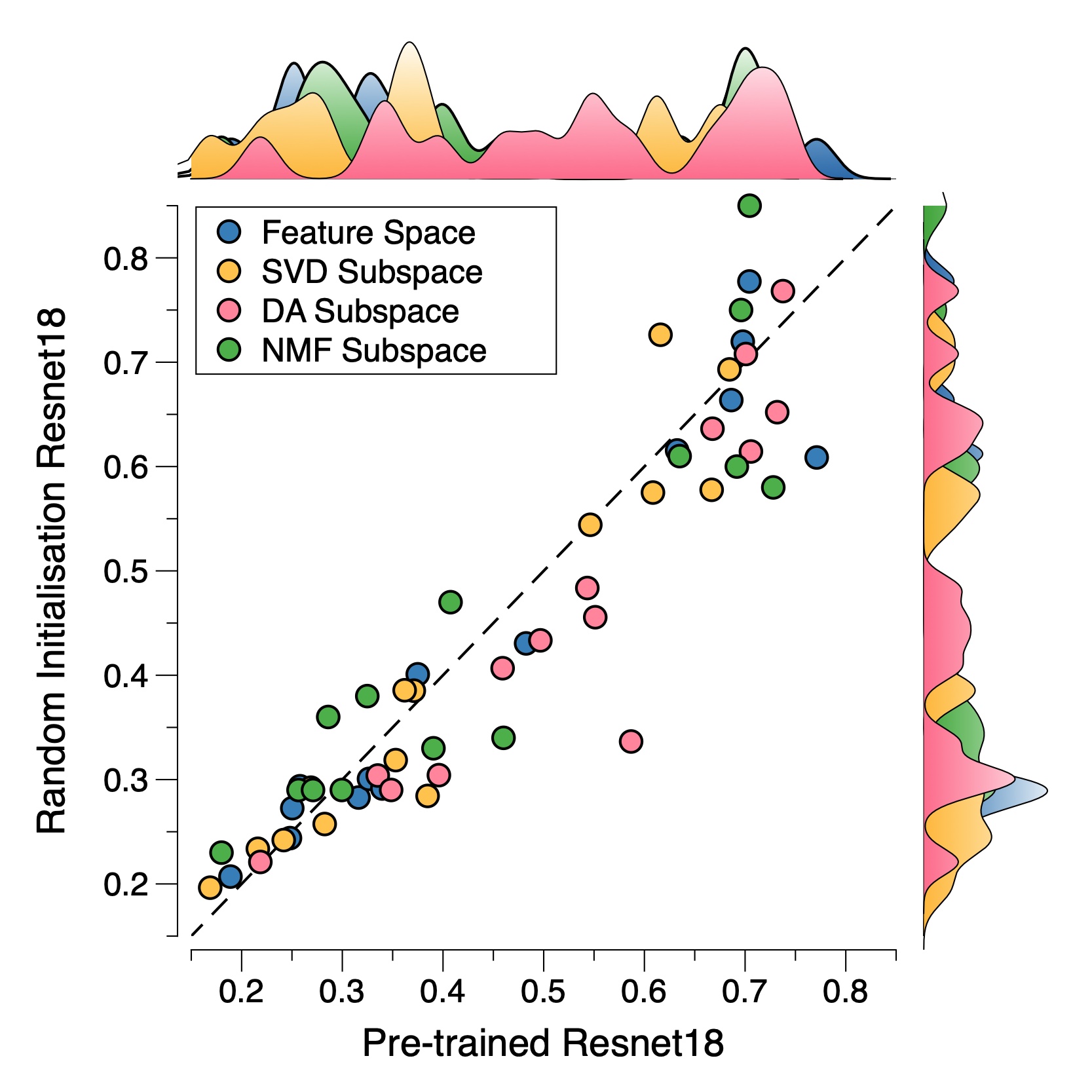}
\caption{Comparison of the use of features derived from pre-trained models against models with random initializations in the few-shot learning framework on the 14 datasets. Information extracted from the pre-trained source models helps in downstream medical tasks, although the fixed random transformations also retain discriminant information.
}
\label{fig:method}
\end{figure}

In the few-shot learning paradigm considered, the pre-trained source model serves as a feature extractor, mapping the medical images into a high dimensional space. To explore the impact of parameters in the model, we compare the classification accuracy from the related subspaces (i.e., feature space, PCA, DA and NMF) in random initialization and pre-trained models.
Figure \ref{fig:method} shows the performance of the pre-trained model and the average of ten random initialization models on all the $14$ datasets. ResNet18 is used as the base feature extractor with various parameters in this experiment. 
As we expected, the features extracted by the pre-trained model retain the good discriminant properties. Surprisingly, the performance of the features extracted by the randomly initialized model and the corresponding subspaces is not significantly degraded, indicating that the same discriminative properties are properly preserved in its extracted features. The DA results in the figure further illustrate this point and prove that subspace perspective provides directions for solving the few-shot learning on medical imaging.

\subsection{Boruta Subspace}\label{Boruta Subspace}
To investigate the performance of feature selection techniques in the few-shot learning framework, we below compare the subspace extracted from the Boruta feature selection method with the dimensionality reduction methods (i.e., SVD, DA and NMF). We follow the Boruta method and extract the related features on the $14$ medical datasets (see results in Appendix \ref{app-boruta}). Figure \ref{fig:boruta-result} presents the classification results comparing the Boruta feature selection method against DA and NMF. It shows that feature selection, like the Boruta method which only selects a subset from the 512 dimensions based on the voting results of a wrapper algorithm around a random forest, generally is not a good choice for the few-shot learning architecture we present. 
Instead of selecting features randomly like Boruta, we prefer to conditionally maintain the original attributes (e.g. discriminability, sparsity and non-negativity) of the data in the subspace. In addition, in terms of the computation time, DA and NMF is dramatically faster than Boruta (needing a high number of iterations), showing the efficiency of the introduced subspace representations.

\begin{figure}[htp]
    \centering
    \subfigure[Boruta vs DA]{
    \centering
    \includegraphics[width=0.35\linewidth]{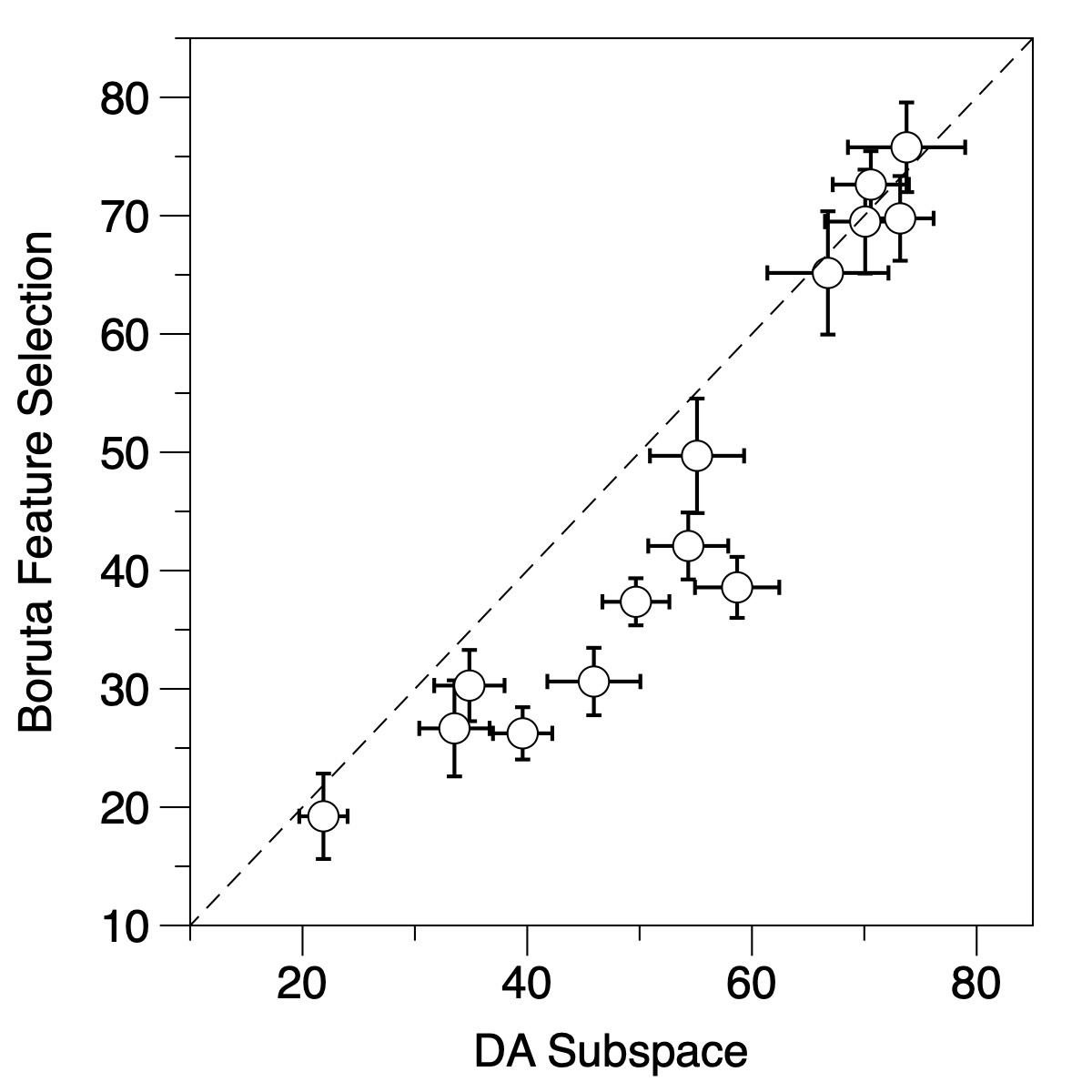}
    }%
    \subfigure[Boruta vs NMF]{
    \centering
    \includegraphics[width=0.35\linewidth]{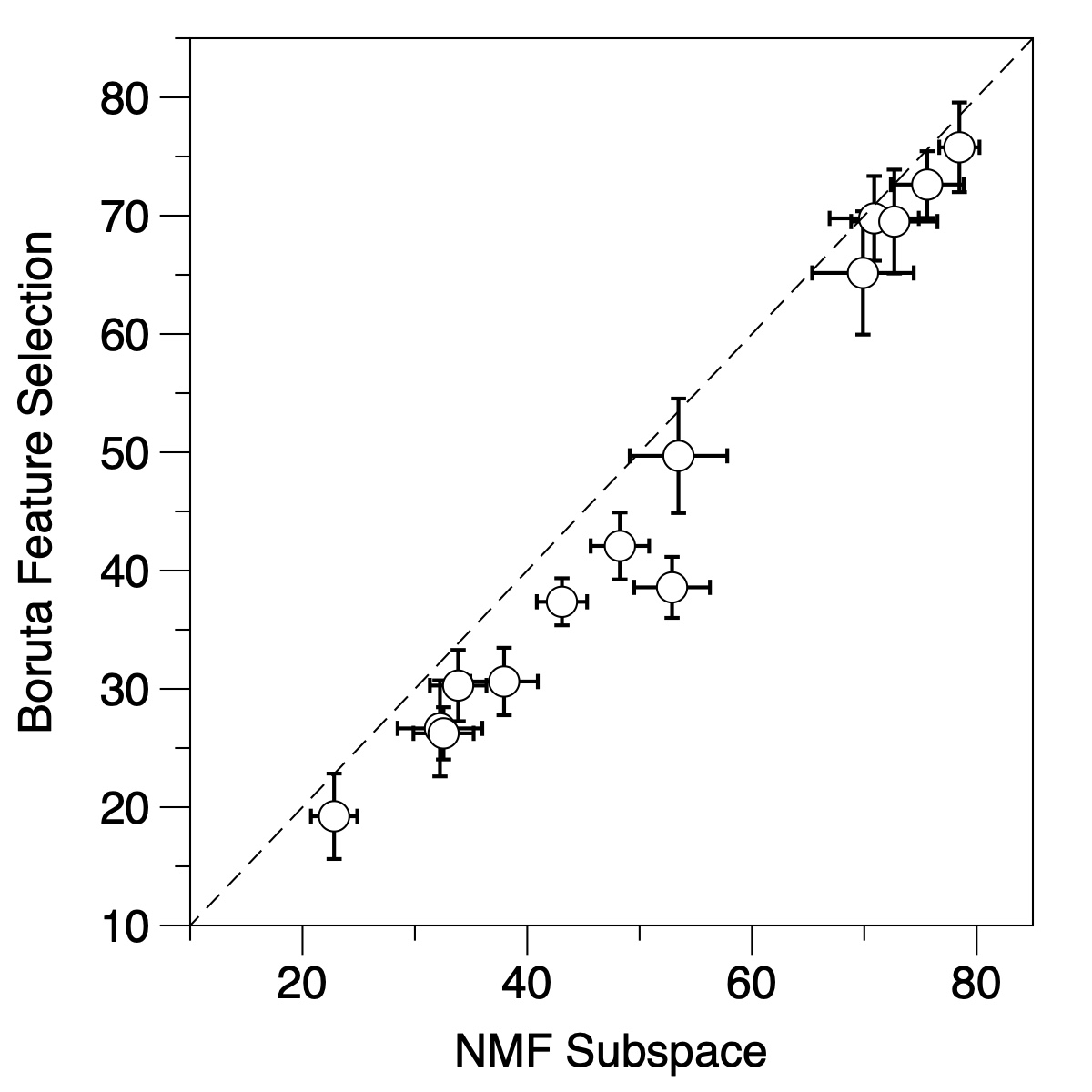}
    }%
\caption{Few-shot learning classification accuracy comparison between Boruta feature selection and the approaches using DA (left panel) and NMF (right panel) subspaces on the 14 distinct medical datasets.}
\label{fig:boruta-result}
\end{figure}

\section{Discussion}\label{Sec:dis}
For few-shot learning with only hundreds of images, comparable in order of magnitude to the feature dimensions (typically $512$ or $1024$ of popular models), dimensionality reduction is essential. While popular method of dimensionality reduction is PCA/SVD, its limitations as a variance preserving approximation suitable for uni-modal data need to be considered. We have addressed this by exploring DA and NMF as alternatives to SVD for few-shot learning in medical imaging.

By presenting the results in the experiment section, we discovered that the subspace obtained by DA is more useful for classification problems than the variance-preserving dimensionality reduction PCA/SVD. DA performs well on multiple disease datasets and effectively distinguishes the classes of disease in the low-dimensional space. However, DA also has some limitations, e.g. the maximum dimension of its subspace is one less than the number of classes for multiclass problems. This limitation is related to the rank of the covariance generated by the dataset. Moreover, DA may not perform ideally with classification when the data information depends on variance rather than the mean.

We also restricted our work on SNMF (supervised NMF) to binary classification problems for which the derivation is readily available. While for multiclass problems, more attempts will be necessary. This is mainly due to the fact that NMF is an inherently unsupervised matrix factorization algorithm and how to properly combine label signals and generate discriminate subspaces remain to be discussed. These, however, do not limit the scope of the conclusions we  reach regarding the desirability of alternatives to the widely used SVD. Future work could be focusing on deriving the solutions to these cases. Additionally, it is also interesting to explore automatic rank selection using information theoretic concepts such as minimum description length considered in \citep{squires2017rank}.

The comparison between feature selection techniques e.g. \citep{li2015robust,tang2020cart} and the dimensionality reduction (i.e., SVD, DA and NMF) reveals that just selecting some specific features is less effective than eliminating less relevant information via dimensionality reduction. Morevover, plain feature selection can be quite unstable and may also be time-consuming. In comparison, since our few-shot learning architecture uses a pre-trained network for feature extraction, it is quite efficient. Most of the time consumed by our few-shot learning architecture is the dimensionality reduction and classification with a simple classifier. Benefiting from the dimensionality reduction, the final classification step is also quite economical.

Finally, it is worth mentioning that in clinical settings the validation and accuracy evaluation of the developed technique in medical imaging are extremely challenging (which is also true for all the related techniques). This is far beyond the lack of data challenge since clinical settings may require the involvement of clinicians, hospitals, patients and even the government, which are all difficult to reach out for individual academics or research groups. Collective effort from all interests is essential to validate/evaluate the practical use of any new method in medical imaging.

\section{Conclusion}\label{Sec:con}
In this paper, we explored two different subspace representations -- DA and NMF -- of features learned from deep neural networks pre-trained on large computer vision datasets, adopted for few-shot learning on small medical imaging datasets. 
Our empirical work is carried out on $14$ different datasets spanning $11$ distinct diseases and four image acquisition modalities. Across these, we demonstrate the following: I) there is a consistent performance advantage on dimensionality reduction in the few-shot learning on medical imaging; II) working with DA derived subspaces gives significant performance gains over PCA/SVD based variance preserving dimensionality reductions, and even when taken at very low dimensions, these gains are statistically significant; and III) NMF-based representation, including its supervised variation, is a viable alternative to SVD-based low dimensional subspaces. NMF also shows a comparable advantage on part-based representation in moderate low dimensions. Overall, the developed few-shot learning framework with the newly introduced subspace representations is a very powerful approach in tackling medical imaging multiclass classification problems.
One of important future avenues could be extending the developed approaches in this work in other fields.

\section*{Acknowledgements}
MN's contribution to this work was funded by Grant EP/S000356/1, Artificial and Augmented Intelligence for Automated Scientific Discovery, Engineering and Physical Sciences Research Council (EPSRC), UK.

\appendices

\begin{figure*}[th]
  \centering
  \includegraphics[width=6.5in]{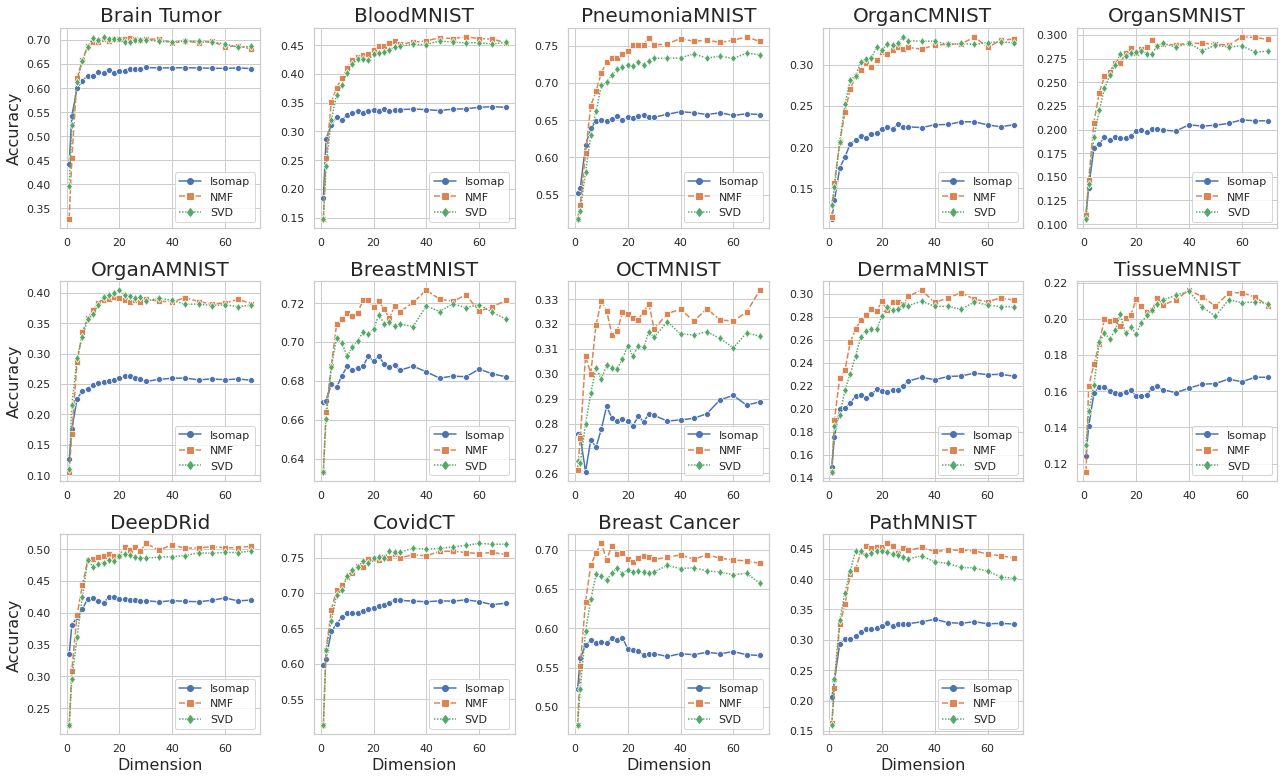}
  \caption{Classification accuracy comparison between the SVD, Isomap and NMF subspaces on 14 datasets with subspace dimensions ranging from 1 to 70. Ten random partitions of the training-test set on each of the 14 datasets are conducted. It shows that in many cases NMF subspace  exhibits slightly better performance in different dimensions than the SVD subspace. Moreover, both SVD and NMF outperform the Isomap subspace by a large margin. }
  \label{all_results}
\end{figure*}

\section{Supervised NMF}\label{Appendix-snmf}
The supervised NMF suggested by 
\citep{leuschner2019supervised}, introducing a logistic regression model into the cost function of NMF in Eq. \eqref{eqn:nmf-min}.  
Let $\boldsymbol{Z}:=\left[1 \mid \boldsymbol{Y X}^{\top}\right] = (\boldsymbol{z}_1, \boldsymbol{z}_2, \cdots, \boldsymbol{z}_N )^\top \in \mathbb{R}^{N\times (p+1)}$. Note that $\boldsymbol{z}_i \in \mathbb{R}^{p+1}, 1\le i \le N$. Considering the binary classification problem, let $\boldsymbol{u} \in\{0,1\}^{N}$ be the vector representing the labels (i.e., 0 or 1) of the given $N$ samples. 
Let $\boldsymbol{\beta} \in \mathbb{R}^{p+1}$. The total loss function of the supervised NMF is defined as
\begin{align} \label{eqn:s-nmf}
\min _{\boldsymbol{K, X} \geq 0, \boldsymbol{\beta}}  \frac{1}{2}  \|\boldsymbol{Y}-\boldsymbol{K X}\|_{F}^{2}+  \frac{\tilde{\lambda}}{N}\Big(\sum_{i=1}^{N} \log \big(1+\mathrm{exp}(\boldsymbol{z}_i^\top \boldsymbol{\beta})\big)  
- \boldsymbol{u}^{\top} \boldsymbol{Z} \boldsymbol{\beta} \Big),
\end{align}
where $\tilde{\lambda} \geq 0$ is the regularization parameter. Problem \eqref{eqn:s-nmf} can be minimised corresponding to $\boldsymbol{K, X}$ and $\boldsymbol{\beta}$ alternatively. Minimisation with respect to $\boldsymbol{K}$ has the same update rules with \eqref{eqn:nmf-update}. Because of the additional included logistic regression term in the loss function, we follow the stochastic gradient descent ADADELTA used in \citep{leuschner2019supervised} to update $\boldsymbol{X}$ and $\boldsymbol{\beta}$. Method ADADELTA leverages the knowledge from previously computed gradients to determine the required gradients $\Delta \boldsymbol{X}$ and $\Delta \boldsymbol{\beta}$ in the next step. 
To guarantee non-negativity of $\boldsymbol{X}$, one projection step, denoted by $\operatorname{proj}(\boldsymbol{X})$, is added to replace all entries of the $\boldsymbol{X}$ less than 0 with a small positive number. The whole update rules for supervised NMF is
\begin{equation}
\begin{aligned}
\boldsymbol{K}   \leftarrow \boldsymbol{K} \circ \frac{\boldsymbol{Y X}^{\top}}{\boldsymbol{K X X^{\top}}}, \ \  
\boldsymbol{X}  \leftarrow \operatorname{proj}(\boldsymbol{X})+\boldsymbol{\Delta X},  \ \
\boldsymbol{\beta}  \leftarrow \boldsymbol{\beta}+\boldsymbol{\Delta \beta}.
\end{aligned}
\end{equation}
For more details please see \citep{leuschner2019supervised}.

\section{Experiment Supplement}
\subsection{Comparison with Manifold Learning}\label{Appendix-nmf-result}
\begin{figure}[th]
\centering
\includegraphics[trim={{.0\linewidth} {0.0\linewidth} {0.0\linewidth} {.0\linewidth}}, clip, width=6.4in]{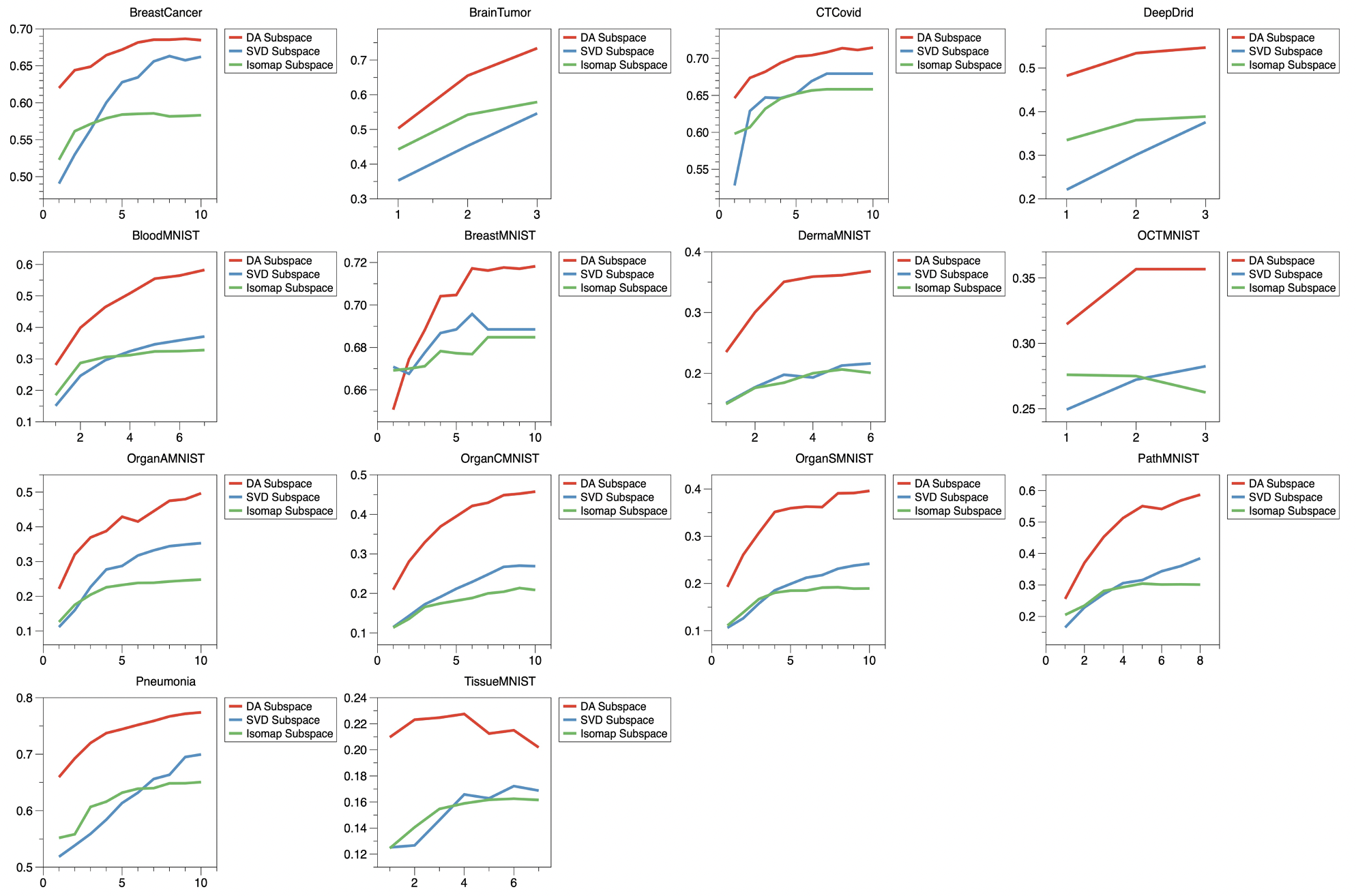}
\caption{Classification accuracy comparison between the SVD, Isomap and DA subspaces on the 14 datasets. 
It shows that DA outperforms both SVD and Isomap  by a large margin.}
\label{isomap}
\end{figure}

A manifold in mathematics is a topological space that resembles Euclidean space. One of the most commonly used manifold learning methods is the nonlinear dimensionality reduction technique Isomap \citep{balasubramanian2002isomap}. Isomap is used for computing a quasi-isometric and low-dimensional embedding of a set of high-dimensional data points.
Comparison between NMF, SVD and Isomap subspaces in terms of classification accuracy on the 14 datasets as the subspace dimensions ranging from 1 to 70 is given in Figure \ref{all_results}.  Ten random partitions of the training-test set on each of the 14 datasets are conducted. It shows that in many cases NMF subspace  exhibits slightly better performance in different dimensions than the SVD subspace. Moreover, both SVD and NMF outperform the Isomap subspace by a large margin.
The comparison between DA, SVD and Isomap subspaces in terms of classification accuracy on all the datasets is given in Figure \ref{isomap}. It shows that Isomap and SVD produce similar results in low dimensions, while DA is able to achieve better results than both Isomap and SVD by a large margin. On the whole, the above results again suggest the viable alternatives of using NMF and DA to the SVD-based low dimensional subspaces.

\subsection{Boruta Results}\label{app-boruta}
Table \ref{table:fea-sel} gives the classification results based on Boruta feature selection and the average number of selected features across ten runs for each of the 14 medical datasets, from which we see that the
Boruta approach is highly unstable (i.e., yielding large deviation).

\begin{table}[!t]
\caption{Few-shot learning classification accuracy by Boruta feature selection on the 14 medical datasets.}
\label{table:fea-sel}
\centering
\setlength{\tabcolsep}{0.8mm}
{\begin{tabular}{c||c|c|c}
\toprule
Datasets & \shortstack{Accuracy} & \shortstack{Selected \\ Features} & Classes
\\
\midrule \midrule
{\tt BreastCancer}\citep{breastcancer} 
&65.16$\pm$5.21&32&2
\\
{\tt BrainTumor}\citep{cheng_2017}
&69.77$\pm$3.58&284&4
\\
{\tt CovidCT}\citep{he2020sample}
&72.63$\pm$2.82&58&2
\\
{\tt DeepDRiD} \citep{DeepDRi}
&49.70$\pm$4.84&106&5
\\
{\tt BloodMNIST}\citep{acevedo2020dataset}
&42.08$\pm$2.83&183&8
\\
{\tt BreastMNIST}\citep{al2020dataset}
&69.50$\pm$4.39&13&2
\\
{\tt DermaMNIST}\citep{tschandl2018ham10000}
&26.66$\pm$4.06&41&7
\\
{\tt OCTMNIST}\citep{kermany2018identifying}
&30.28$\pm$3.01&19&4
\\
{\tt OrganAMNIST}\citep{bilic2019liver}
&37.36$\pm$1.99&170&11
\\
{\tt OrganCMNIST}\citep{bilic2019liver}
&30.62$\pm$2.85&125&11
\\
{\tt OrganSMNIST}\citep{bilic2019liver}
&26.24$\pm$2.21&100&11
\\
{\tt PathMNIST}\citep{kather2019predicting}
&38.58$\pm$2.58&197&9
\\
{\tt PneumoniaMNIST}\citep{kermany2018identifying}
&75.78$\pm$3.79&71&2
\\
{\tt TissueMNIST}\citep{woloshuk2021situ}
&19.23$\pm$3.61&26&8\\
\midrule
\end{tabular}}
\end{table}

\begin{table}[]
\setlength{\tabcolsep}{0.7mm}{
\caption{Classification accuracy comparison between the prototypical network and the few-shot learning with subspace feature representations. In particular, $5$ samples from each class in each dataset are used for training, i.e., forming the ``$C$-way 5-shot" setting (recall $C$ is the number of classes in each dataset). Dim stands for dimensions. }
\label{fsw-all}
\resizebox{\textwidth}{!}{
\begin{tabular}{ccccccccccc}
\hline
\multicolumn{11}{c}{C-way 5-shot  Accuracy(\%)} \\ \hline
\multirow{3}{*}{\begin{tabular}[c]{@{}c@{}}Data\end{tabular}} &
  \multicolumn{10}{c}{Methods} \\ \cline{2-11} 
 &
  \multirow{2}{*}{\begin{tabular}[c]{@{}c@{}}Prototypical\\ Network\end{tabular}} &
  \multicolumn{9}{c}{Few-shot Learning with Subspace Feature Representations (\textbf{Ours})} \\ \cline{3-11} 
 &
   &
  Feature Space &
  \multicolumn{1}{l}{Subspaces} &
  2 Dim &
  5 Dim &
  10 Dim &
  20 Dim &
  30 Dim &
  40 Dim &
  50 Dim \\ \hline
\multirow{3}{*}{\begin{tabular}[c]{@{}c@{}}CovidCT\\ (2 classes)\end{tabular}} &
  \multirow{3}{*}{53.33$\pm$6.99} &
  \multirow{3}{*}{52.22$\pm$5.88} &
  SVD &
  
  49.78$\pm$10.95 &
  54.56$\pm$5.80 &
  52.11$\pm$5.32 &
  52.11$\pm$5.32 &
  52.11$\pm$5.32 &
  52.11$\pm$5.32 &
  52.11$\pm$5.32 \\
 &
   &
   &
  NMF &
  53.00$\pm$10.05 &
  \textbf{56.89$\pm$8.07} &
  56.44$\pm$4.89 &
  55.89$\pm$6.61 &
  55.78$\pm$6.36 &
  56.33$\pm$6.09 &
  56.44$\pm$6.08 \\
  \cline{4-11}
   &
   &
   &
  DA &
  \multicolumn{7}{c}{51.91$\pm$3.41 (10 Dim)}
   \\
  \hline
\multirow{3}{*}{\begin{tabular}[c]{@{}c@{}}BreastCancer\\ (2 classes)\end{tabular}} &
  \multirow{3}{*}{72.33$\pm$8.68} &
  \multirow{3}{*}{70.22$\pm$8.68} &
  SVD &
  66.89$\pm$5.83 &
  70.33$\pm$8.61 &
  70.11$\pm$9.07 &
  70.11$\pm$9.07 &
  70.11$\pm$9.07 &
  70.11$\pm$9.07 &
  70.11$\pm$9.07 \\
 &
   &
   &
  NMF &
  68.78$\pm$7.83 &
  72.44$\pm$7.71 &
  72.78$\pm$9.69 &
  72.11$\pm$7.82 &
  \textbf{73.00$\pm$8.93} &
  71.33$\pm$9.27 &
  72.44$\pm$8.37 \\ 
  \cline{4-11}
     &
   &
   &
  DA &
  \multicolumn{7}{c}{62.69$\pm$4.43  (10 Dim)}
   \\
 
  \hline
\multirow{3}{*}{\begin{tabular}[c]{@{}c@{}}PneumoniaMNIST\\ (2 classes)\end{tabular}} &
  \multirow{3}{*}{64.33$\pm$8.17} &
  \multirow{3}{*}{66.56$\pm$9.79} &
  SVD &
  61.56$\pm$9.34 &
  66.11$\pm$7.32 &
  66.56$\pm$8.92 &
  66.56$\pm$8.92 &
  66.56$\pm$8.92 &
  66.56$\pm$8.92 &
  66.56$\pm$8.92 \\
 &
   &
   &
  NMF &
  62.67$\pm$9.44 &
  70.33$\pm$6.38 &
  72.00$\pm$6.46 &
  72.89$\pm$6.03 &
  72.56$\pm$5.75 &
  73.33$\pm$7.42 &
  \textbf{73.78$\pm$6.78} \\ 
   \cline{4-11}
   &
   &
   &
   DA &
   \multicolumn{7}{c}{68.45$\pm$5.79 (10 Dim)}
   \\
   \hline
\multirow{3}{*}{\begin{tabular}[c]{@{}c@{}}BreastMNIST\\ (2 classes)\end{tabular}} &
  \multirow{3}{*}{54.00$\pm$9.40} &
  \multirow{3}{*}{54.11$\pm$6.43} &
  SVD &
  54.00$\pm$7.14 &
  54.11$\pm$6.67 &
  54.67$\pm$7.47 &
  54.67$\pm$7.47 &
  54.67$\pm$7.47 &
  54.67$\pm$7.47 &
  54.67$\pm$7.47 \\
 &
   &
   &
  NMF &
  53.00$\pm$7.70 &
  59.11$\pm$4.63 &
  60.33$\pm$5.17 &
  60.56$\pm$3.66 &
  59.00$\pm$6.84 &
  60.00$\pm$6.19 &
  \textbf{62.33$\pm$5.45} \\
  \cline{4-11}
  &
  &
  &
  DA
  &
  \multicolumn{7}{c}{59.98$\pm$8.08 (Dim = 10)}
  \\ 
  \hline
\multirow{3}{*}{\begin{tabular}[c]{@{}c@{}}DeepDRid\\ (5 classes)\end{tabular}} &
  \multirow{3}{*}{30.07$\pm$7.87} &
  \multirow{3}{*}{29.47$\pm$4.91} &
  SVD &
  27.58$\pm$4.53 &
  28.03$\pm$5.59 &
  29.42$\pm$4.75 &
  29.28$\pm$4.95 &
  29.47$\pm$4.98 &
  29.47$\pm$4.98 &
  29.47$\pm$4.98 \\
 &
   &
   &
  NMF &
  28.92$\pm$4.51 &
  30.03$\pm$5.59 &
  31.07$\pm$4.48 &
  31.27$\pm$4.59 &
  31.02$\pm$4.54 &
  \textbf{31.47$\pm$4.17} &
  31.37$\pm$4.07 \\
  \cline{4-11}
     &
   &
   &
  DA &
   \multicolumn{7}{c}{\textbf{31.63$\pm$7.07}  (4 Dim)}
   \\
   \hline
\multirow{3}{*}{\begin{tabular}[c]{@{}c@{}}BrainTumor\\ (4 classes)\end{tabular}} &
  \multirow{3}{*}{34.67$\pm$4.82} &
  \multirow{3}{*}{35.42$\pm$5.29} &
  SVD &
  33.12$\pm$4.69 &
  34.54$\pm$4.84 &
  35.00$\pm$5.29 &
  35.42$\pm$5.23 &
  35.42$\pm$5.23 &
  35.42$\pm$5.23 &
  35.42$\pm$5.23 \\
 &
   &
   &
  NMF &
  33.96$\pm$5.45 &
  36.04$\pm$4.92 &
  36.88$\pm$4.73 &
  37.33$\pm$5.20 &
  \textbf{37.87$\pm$5.10} &
  37.50$\pm$5.10 &
  37.46$\pm$4.87 \\ 
  \cline{4-11}
     &
   &
   &
  DA &
   \multicolumn{7}{c}{\textbf{61.88$\pm$5.50} (3 Dim)}
   \\
   \hline
\multirow{3}{*}{\begin{tabular}[c]{@{}c@{}}BloodMNIST\\ (8 classes)\end{tabular}} &
  \multirow{3}{*}{47.58$\pm$5.24} &
  \multirow{3}{*}{48.29$\pm$2.52} &
  SVD &
  36.29$\pm$2.98 &
  46.42$\pm$4.65 &
  47.62$\pm$2.49 &
  48.06$\pm$3.15 &
  \textbf{48.33$\pm$2.45} &
  48.33$\pm$2.63 &
  48.33$\pm$2.63 \\
 &
   &
   &
  NMF &
  36.60$\pm$3.77 &
  47.33$\pm$4.96 &
  46.94$\pm$4.42 &
  47.02$\pm$3.87 &
  46.23$\pm$2.71 &
  45.90$\pm$3.65 &
  46.58$\pm$3.56 \\ 
  \cline{4-11}
     &
   &
   &
  DA &
  \multicolumn{7}{c}{\textbf{54.33$\pm$8.08} (7 Dim)}
   \\
   \hline
\multirow{3}{*}{\begin{tabular}[c]{@{}c@{}}DermaMNIST\\ (7 classes)\end{tabular}} &
  \multirow{3}{*}{26.57$\pm$4.69} &
  \multirow{3}{*}{25.83$\pm$3.57} &
  SVD &
  21.43$\pm$4.61 &
  25.07$\pm$3.63 &
  25.19$\pm$3.28 &
  25.33$\pm$3.67 &
  25.71$\pm$3.57 &
  25.76$\pm$3.58 &
  25.76$\pm$3.58 \\
 &
   &
   &
  NMF &
  22.71$\pm$3.49 &
  27.45$\pm$4.02 &
  26.19$\pm$3.36 &
  26.55$\pm$3.16 &
  27.86$\pm$3.40 &
  26.83$\pm$2.92 &
  \textbf{27.79$\pm$2.98} \\ 
  \cline{4-11}
  &
  &
  &
  DA &
   \multicolumn{7}{c}{\textbf{31.14$\pm$5.54} (6 Dim)}
   \\
   \hline
\multirow{3}{*}{\begin{tabular}[c]{@{}c@{}}OCTMNIST\\ (4 classes)\end{tabular}} &
  \multirow{3}{*}{28.64$\pm$3.95} &
  \multirow{3}{*}{29.79$\pm$3.02} &
  SVD &
  26.08$\pm$3.80 &
  28.50$\pm$4.31 &
  29.42$\pm$2.94 &
  29.83$\pm$3.43 &
  29.83$\pm$3.43 &
  29.83$\pm$3.43 &
  29.83$\pm$3.43 \\
 &
   &
   &
  NMF &
  26.12$\pm$3.63 &
  31.17$\pm$2.52 &
  \textbf{32.17$\pm$3.12} &
  31.75$\pm$3.20 &
  32.13$\pm$3.49 &
  31.04$\pm$3.16 &
  32.12$\pm$3.57 \\ 
  \cline{4-11} 
   &
   &
   &
  DA & \multicolumn{7}{c}{\textbf{32.92$\pm$6.68} ($3$ Dim)}
   \\
   \hline
\multirow{3}{*}{\begin{tabular}[c]{@{}c@{}}OrganAMNIST\\ (11 classes)\end{tabular}} &
  \multirow{3}{*}{47.94$\pm$2.39} &
  \multirow{3}{*}{52.50$\pm$3.03} &
  SVD &
  34.26$\pm$3.28 &
  45.52$\pm$3.92 &
  50.03$\pm$3.36 &
  52.73$\pm$3.29 &
  52.74$\pm$3.48 &
  52.89$\pm$3.15 &
  52.71$\pm$3.22 \\
 &
   &
   &
  NMF &
  33.70$\pm$3.38 &
  45.39$\pm$2.67 &
  53.00$\pm$3.31 &
  53.88$\pm$3.38 &
  \textbf{53.88$\pm$3.26} &
  52.97$\pm$2.82 &
  53.58$\pm$3.99 \\ 
  \cline{4-11}
   &
   &
   &
  DA &
   \multicolumn{7}{c}{\textbf{60.94$\pm$3.86} (10 Dim)}
   \\
   \hline
\multirow{3}{*}{\begin{tabular}[c]{@{}c@{}}OrganCMNIST\\ (11 classes)\end{tabular}} &
  \multirow{3}{*}{48.55$\pm$4.11} &
  \multirow{3}{*}{50.23$\pm$3.84} &
  SVD &
  30.61$\pm$3.48 &
  42.12$\pm$4.15 &
  48.68$\pm$3.71 &
  49.70$\pm$4.23 &
  50.18$\pm$4.46 &
  50.18$\pm$4.02 &
  50.08$\pm$4.25 \\
 &
   &
   &
  NMF &
  29.86$\pm$2.77 &
  40.55$\pm$3.83 &
  49.55$\pm$3.67 &
  51.92$\pm$4.72 &
  \textbf{52.29$\pm$4.40} &
  51.26$\pm$4.55 &
  50.91$\pm$3.55 \\ 
  \cline{4-11}
     &
   &
   &
  DA &
   \multicolumn{7}{c}{\textbf{60.62$\pm$3.10} (10 Dim)}
   \\
   \hline
\multirow{3}{*}{\begin{tabular}[c]{@{}c@{}}OrganSMNIST\\ (11 classes)\end{tabular}} &
  \multirow{3}{*}{34.67$\pm$4.21} &
  \multirow{3}{*}{36.95$\pm$3.69} &
  SVD &
  25.74$\pm$2.64 &
  33.79$\pm$3.32 &
  35.70$\pm$4.01 &
  36.52$\pm$3.86 &
  \textbf{37.18$\pm$3.43} &
  37.05$\pm$3.45 &
  36.94$\pm$3.58 \\
 &
   &
   &
  NMF &
  24.67$\pm$2.41 &
  33.97$\pm$3.04 &
  36.06$\pm$3.24 &
  36.45$\pm$3.11 &
  35.48$\pm$2.97 &
  35.53$\pm$3.34 &
  34.67$\pm$3.80 \\ 
  \cline{4-11}
     &
   &
   &
  DA &
   \multicolumn{7}{c}{\textbf{41.23$\pm$2.62} (10 Dim)}
   \\
   \hline
\multirow{3}{*}{\begin{tabular}[c]{@{}c@{}}PathMINST\\ (9 classes)\end{tabular}} &
  \multirow{3}{*}{36.22$\pm$4.77} &
  \multirow{3}{*}{37.69$\pm$3.53} &
  SVD &
  28.83$\pm$5.01 &
  36.17$\pm$4.08 &
  37.72$\pm$3.45 &
  37.96$\pm$3.72 &
  37.50$\pm$3.53 &
  37.76$\pm$3.64 &
  37.74$\pm$3.59 \\
 &
   &
   &
  NMF &
  28.67$\pm$4.26 &
  36.69$\pm$3.04 &
  \textbf{39.76$\pm$3.10} &
  39.07$\pm$3.33 &
  38.67$\pm$2.42 &
  37.54$\pm$2.95 &
  37.04$\pm$3.29 \\ 
  \cline{4-11} 
   &
   &
   &
  DA 
  &
  \multicolumn{7}{c}{\textbf{41.43$\pm$4.19} ($8$ Dim)}
   \\
   \hline
\multirow{3}{*}{\begin{tabular}[c]{@{}c@{}}TissueMNIST\\ (8 classes)\end{tabular}} &
  \multirow{3}{*}{24.42$\pm$3.67} &
  \multirow{3}{*}{23.65$\pm$2.22} &
  SVD &
  19.02$\pm$3.15 &
  20.94$\pm$2.74 &
  22.60$\pm$2.39 &
  23.21$\pm$2.53 &
  23.83$\pm$2.44 &
  23.75$\pm$2.21 &
  23.75$\pm$2.21 \\
 &
   &
   &
  NMF &
  18.50$\pm$3.63 &
  22.40$\pm$2.43 &
  \textbf{24.58$\pm$2.22} &
  24.23$\pm$2.39 &
  23.31$\pm$2.58 &
  22.96$\pm$2.63 &
  23.02$\pm$1.72 \\ 
  \cline{4-11}
     &
   &
   &
  DA &
   \multicolumn{7}{c}{\textbf{38.35$\pm$5.59} (7 Dim)}
   \\\hline
\end{tabular}}
}
\end{table}

\subsection{Classification Results via SVM}\label{svm-nmf-da}
\begin{figure*}[!t]
\centering
\includegraphics[trim={{.0\linewidth} {.0\linewidth} {.0\linewidth} {.00\linewidth}}, clip, width=4.0in, height=1.6in]{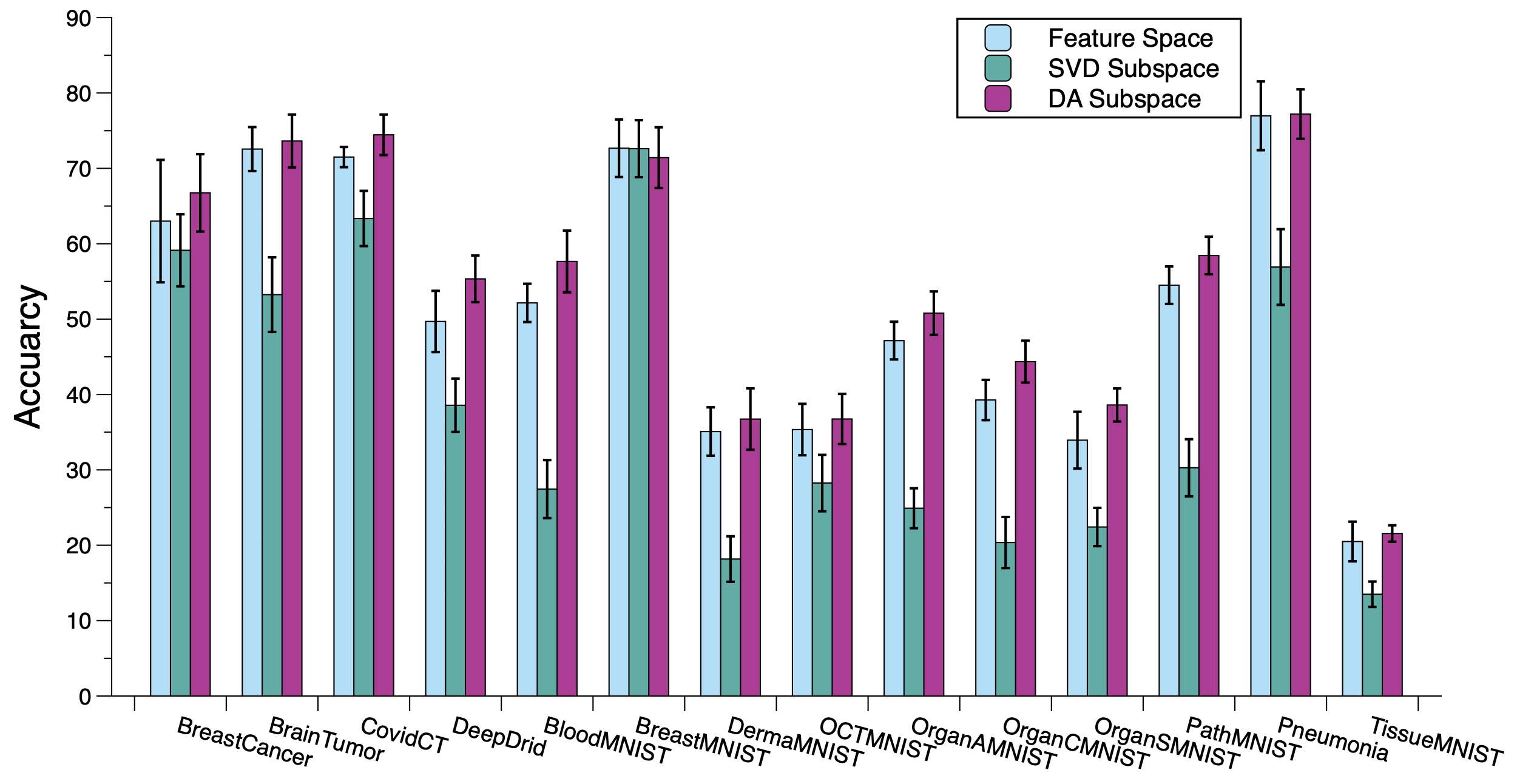}
\includegraphics[width=2.1in,height=2.0in]{NMF/nmf.jpg}
\caption{Few-shot learning performance using SVM as classifier on the 14 medical datasets. Left panel: results obtained by using the original feature space, the SVD subspace and the DA subspace. Right panel: results obtained by using the features in 30-dimensional subspace derived by SVD and NMF.} 
\label{fig:SVM-da-nmf}
\end{figure*}

Figure \ref{fig:SVM-da-nmf} depicts the performance of applying SVM as a classifier in the developed few-shot learning framework  using the same setup as KNN on the 14 medical datasets. In particular, we implemented both the rbf kernel and the linear kernel for SVM.
The C and gamma parameters were also fine-tuned for different datasets. The same dimensions as NMF and DA were applied to SVD to ensure fair comparison. We discovered that the rbf kernel worked better in comparisons between SVD and DA in low dimensions, and the linear kernel performed better in comparisons between SVD and NMF in medium dimensions. On the whole, consistent results were obtained by using SVM and KNN as classifiers on the 14 distinct medical datasets, indicating that the developed few-shot learning architecture is robust to the choice of classifiers.

\subsection{T-SNE Visualization}\label{TSNE visualisation}

\begin{figure}[htp]
    \centering
    \subfigure[Original feature space]{
    \centering
    \includegraphics[width=0.35\linewidth]{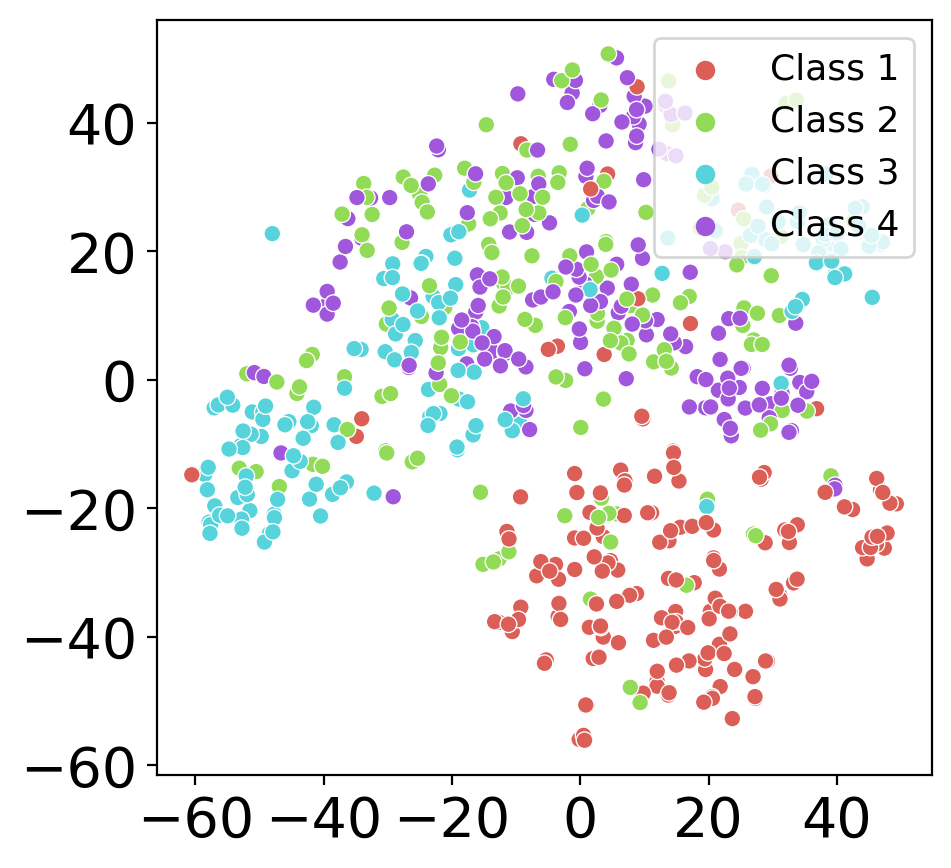}
    \label{tsne_ori}
    }%
    \subfigure[DA]{
    \centering
    \includegraphics[width=0.35\linewidth]{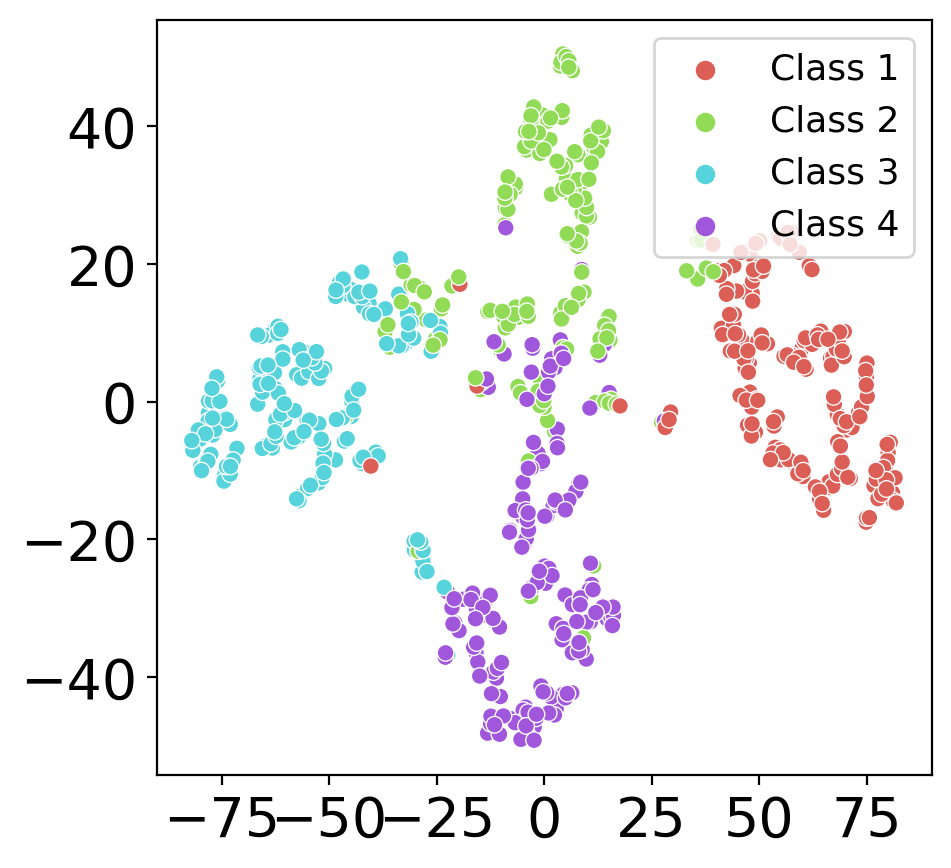}
    \label{tsne_da}
    }%
    
    \subfigure[SVD]{
    \centering
    \includegraphics[width=0.35\linewidth]{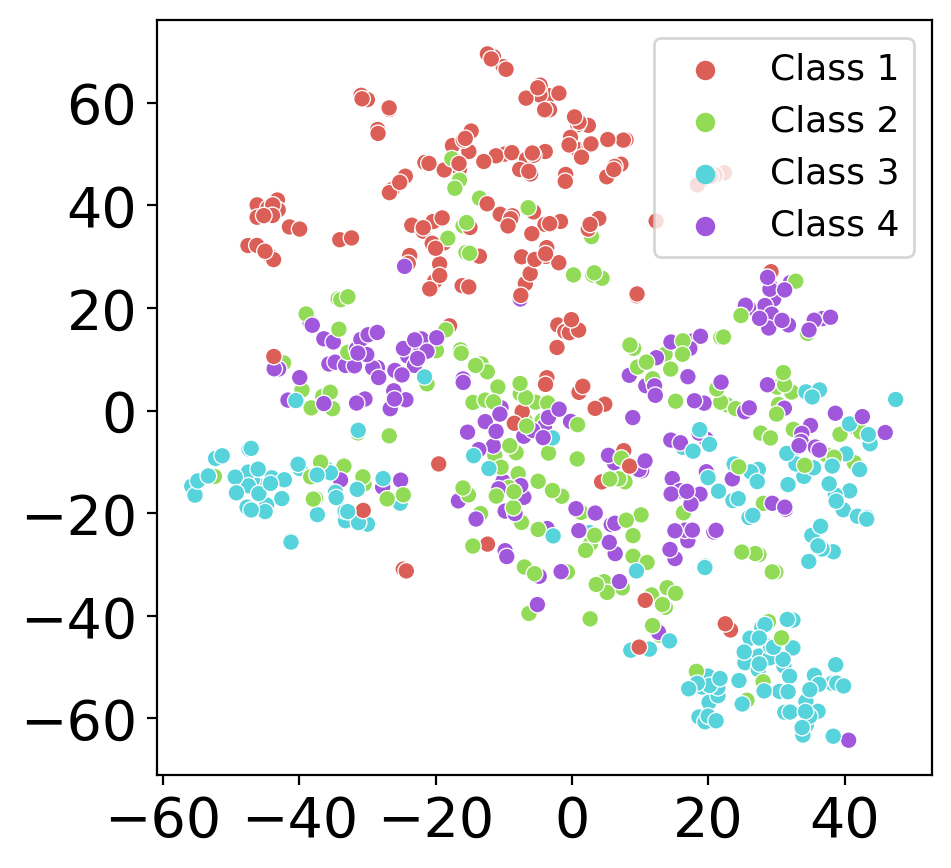}
    \label{tsne_svd}
    }%
    \centering
    \subfigure[NMF]{
    \centering
    \includegraphics[width=0.35\linewidth]{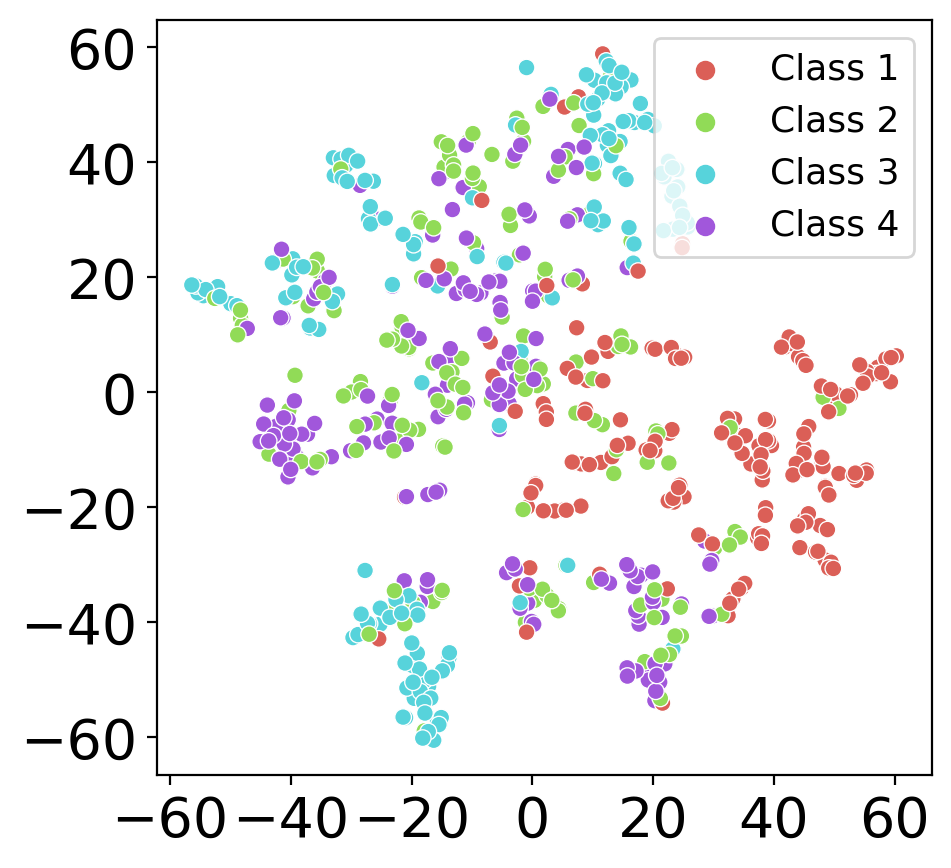}
    \label{tsne_nmf}
    }%
    \centering
    \caption{Subspace visualization by T-SNE on the brain tumour dataset. (a)--(d): the results regarding the originally feature space, DA subspace, 30-dimensional subspace derived by SVD and NMF, respectively.} 
    \label{tsne-all}
\end{figure}

For better visual validation,
Figure \ref{tsne-all} shows an example of visualising the subspaces utilising the T-SNE visualization technique (built-in function in Python) on the brain tumour dataset. Figure \ref{tsne_ori} and Figure \ref{tsne_da} are the 2D feature visualization of the original feature space (i.e., features extracted by the pre-trained network) and the DA subspace, respectively. Figure \ref{tsne_svd} and Figure \ref{tsne_nmf} show the features that are projected to the 30-dimensional subspace by SVD and NMF, respectively.  These plots visually prove that the DA and NMF subspaces are indeed viable alternatives to SVD.

\subsection{Comparison with the Prototypical Network}
\label{Prototypical}

Below we compare our techniques with the well-known few-shot learning algorithm, the prototypical network \citep{snell2017prototypical}, on all the 14 medical datasets, see Table \ref{fsw-all}.
The architecture of the prototypical network used in this experiment is the same as the one in \citep{snell2017prototypical}, which is composed of four convolution blocks and has been trained on the omniglot dataset \citep{lake2019omniglot} via SGD with Adam optimiser \citep{kingma2014adam} and obtained 99\% accuracy in the 5-shot scenario. 
Each block in this network is comprised of a 64-filter $3 \times 3$ convolution, batch normalisation layer, a ReLU nonlinearity and a $2 \times 2$ max-pooling layer. 
The classification accuracy reported in Table \ref{fsw-all} is averaging over 10 randomly generated episodes from the test set.
The setting for the test experiment is "$C$-way $5$-shot," where $5$ samples are given for each class in the support set and $C$ is the number of classes in each dataset. To validate the final performance, 15 query images per class are provided. Recall that the original feature space represents the features extracted by the network without dimensionality reduction.  The dimensions for the NMF/SVD subspaces are chosen as $2$, $5$, $10$, $20$, $30$, $40$ and $50$, separately. The dimensions of the DA subspace are chosen as $10$ for the binary classification problem and $(C-1)$ for the multiclass classification problem. 
The results in Table \ref{fsw-all} demonstrate that our method outperforms all the other methods, i.e., the prototypical network and the ones with the original feature space and SVD.

 \begin{figure}[htp]
\subfigure[{\tt BloodMNIST}]{
\centering
\includegraphics[width=0.48\linewidth]{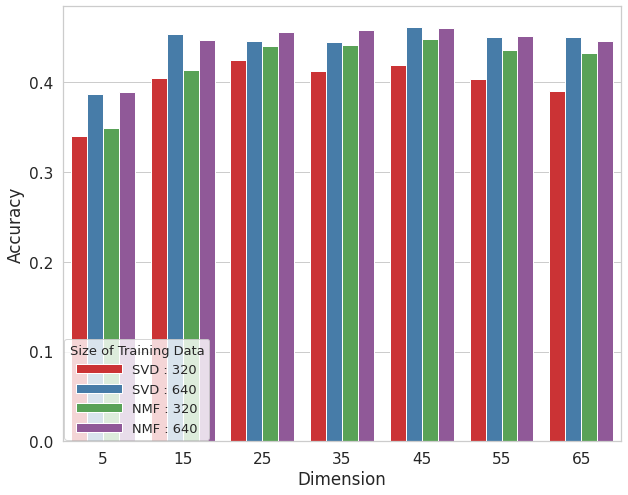}
}%
\subfigure[{\tt DeepDRid}]{
\centering
\includegraphics[width=0.48\linewidth]{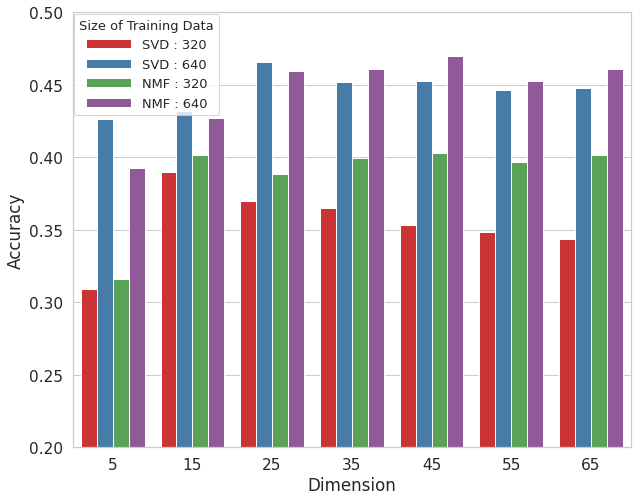}
}%
\caption{Comparison between NMF and SVD subspaces in terms of classification accuracy corresponding to different dataset size as the subspace dimension changes. Datasets  {\tt BloodMNIST} with eight classes and {\tt DeepDRid} with five classes are used in panels (a) and (b), respectively. }
\label{Impact-size-dimensionality-fig}
\end{figure}

\subsection{Impact of the Dataset Size and Dimensionality}
\label{Impact-size-dimensionality}

Figure \ref{Impact-size-dimensionality-fig} shows the impact of the dataset size on the classification accuracy of NMF and SVD as the subspace dimension changes, where datasets {\tt BloodMNIST} with eight classes and {\tt DeepDRid} with five classes are used. The setting in Figure \ref{Impact-size-dimensionality-fig} is the same as the one used in Figure \ref{datavolumn}. It shows that, in the multiclass classification problem, fewer categories will result in higher accuracy, and SVD suffers from dimension changes more in the datasets with small size. 
Without enough data, SVD could not extract specific precise features that match the target categories or may extract insignificant features as the dimension increases, resulting in low classification performance (e.g. see Figure \ref{Impact-size-dimensionality-fig}(b)). In contrast, the results of NMF are robust to dimension changes in datasets with different size. This great performance of NMF benefits from its part-based representation maximising the features preserved in the subspace.

\bibliographystyle{unsrtnat}

\bibliography{main}

\begin{thebibliography}{10}

\bibitem{krizhevsky2012imagenet}
Alex Krizhevsky, Ilya Sutskever, and Geoffrey~E Hinton.
\newblock Imagenet classification with deep convolutional neural networks.
\newblock {\em Advances in neural information processing systems}, 25, 2012.

\bibitem{fukushima1982neocognitron}
Kunihiko Fukushima and Sei Miyake.
\newblock Neocognitron: A self-organizing neural network model for a mechanism
  of visual pattern recognition.
\newblock In {\em Competition and cooperation in neural nets}, pages 267--285.
  Springer, 1982.

\bibitem{pinto2008real}
Nicolas Pinto, David~D Cox, and James~J DiCarlo.
\newblock Why is real-world visual object recognition hard?
\newblock {\em PLoS computational biology}, 4(1):e27, 2008.

\bibitem{GUPTA2021100057}
Abhishek Gupta, Alagan Anpalagan, Ling Guan, and Ahmed~Shaharyar Khwaja.
\newblock Deep learning for object detection and scene perception in
  self-driving cars: Survey, challenges, and open issues.
\newblock {\em Array}, 10:100057, 2021.

\bibitem{zhu2017target}
Yuke Zhu, Roozbeh Mottaghi, Eric Kolve, Joseph~J Lim, Abhinav Gupta,
  Li~Fei-Fei, and Ali Farhadi.
\newblock Target-driven visual navigation in indoor scenes using deep
  reinforcement learning.
\newblock In {\em 2017 IEEE international conference on robotics and automation
  (ICRA)}, pages 3357--3364. IEEE, 2017.

\bibitem{liu2019comparison}
Xiaoxuan Liu, Livia Faes, Aditya~U Kale, Siegfried~K Wagner, Dun~Jack Fu, Alice
  Bruynseels, Thushika Mahendiran, Gabriella Moraes, Mohith Shamdas, Christoph
  Kern, et~al.
\newblock A comparison of deep learning performance against health-care
  professionals in detecting diseases from medical imaging: a systematic review
  and meta-analysis.
\newblock {\em The lancet digital health}, 1(6):e271--e297, 2019.

\bibitem{pesapane2018artificial}
Filippo Pesapane, Marina Codari, and Francesco Sardanelli.
\newblock Artificial intelligence in medical imaging: threat or opportunity?
  radiologists again at the forefront of innovation in medicine.
\newblock {\em European radiology experimental}, 2(1):1--10, 2018.

\bibitem{castro2020causality}
Daniel~C Castro, Ian Walker, and Ben Glocker.
\newblock Causality matters in medical imaging.
\newblock {\em Nature Communications}, 11(1):1--10, 2020.

\bibitem{kendall2017uncertainties}
Alex Kendall and Yarin Gal.
\newblock What uncertainties do we need in bayesian deep learning for computer
  vision?
\newblock {\em Advances in neural information processing systems}, 30, 2017.

\bibitem{LUNDERVOLD2019102}
Alexander~Selvikvåg Lundervold and Arvid Lundervold.
\newblock An overview of deep learning in medical imaging focusing on mri.
\newblock {\em Zeitschrift für Medizinische Physik}, 29(2):102--127, 2019.
\newblock Special Issue: Deep Learning in Medical Physics.

\bibitem{kim2017few}
Mijung Kim, Jasper Zuallaert, and Wesley De~Neve.
\newblock Few-shot learning using a small-sized dataset of high-resolution
  fundus images for glaucoma diagnosis.
\newblock In {\em Proceedings of the 2nd international workshop on multimedia
  for personal health and health care}, pages 89--92, 2017.

\bibitem{althnian2021impact}
Alhanoof Althnian, Duaa AlSaeed, Heyam Al-Baity, Amani Samha, Alanoud~Bin Dris,
  Najla Alzakari, Afnan Abou~Elwafa, and Heba Kurdi.
\newblock Impact of dataset size on classification performance: An empirical
  evaluation in the medical domain.
\newblock {\em Applied Sciences}, 11(2):796, 2021.

\bibitem{raghu2019transfusion}
Maithra Raghu, Chiyuan Zhang, Jon Kleinberg, and Samy Bengio.
\newblock Transfusion: Understanding transfer learning for medical imaging.
\newblock {\em Advances in neural information processing systems}, 32, 2019.

\bibitem{CAROPPO2021101852}
Andrea Caroppo, Alessandro Leone, and Pietro Siciliano.
\newblock Deep transfer learning approaches for bleeding detection in endoscopy
  images.
\newblock {\em Computerized Medical Imaging and Graphics}, 88:101852, 2021.

\bibitem{zhang2021weakly}
Dingwen Zhang, Junwei Han, Gong Cheng, and Ming-Hsuan Yang.
\newblock Weakly supervised object localization and detection: A survey.
\newblock {\em IEEE Transactions on Pattern Analysis and Machine Intelligence},
  44(9):5866--5885, 2021.

\bibitem{huang2021unsupervised}
Hsin-Ping Huang, Krishna~C Puvvada, Ming Sun, and Chao Wang.
\newblock Unsupervised and semi-supervised few-shot acoustic event
  classification.
\newblock In {\em ICASSP 2021-2021 IEEE International Conference on Acoustics,
  Speech and Signal Processing (ICASSP)}, pages 331--335. IEEE, 2021.

\bibitem{peng2019few}
Zhimao Peng, Zechao Li, Junge Zhang, Yan Li, Guo-Jun Qi, and Jinhui Tang.
\newblock Few-shot image recognition with knowledge transfer.
\newblock In {\em Proceedings of the IEEE/CVF International Conference on
  Computer Vision}, pages 441--449, 2019.

\bibitem{tang2020blockmix}
Hao Tang, Zechao Li, Zhimao Peng, and Jinhui Tang.
\newblock Blockmix: meta regularization and self-calibrated inference for
  metric-based meta-learning.
\newblock In {\em Proceedings of the 28th ACM international conference on
  multimedia}, pages 610--618, 2020.

\bibitem{wang2020generalizing}
Yaqing Wang, Quanming Yao, James~T Kwok, and Lionel~M Ni.
\newblock Generalizing from a few examples: A survey on few-shot learning.
\newblock {\em ACM computing surveys (csur)}, 53(3):1--34, 2020.

\bibitem{argueso2020few}
David Arg{\"u}eso, Artzai Picon, Unai Irusta, Alfonso Medela, Miguel~G
  San-Emeterio, Arantza Bereciartua, and Aitor Alvarez-Gila.
\newblock Few-shot learning approach for plant disease classification using
  images taken in the field.
\newblock {\em Computers and Electronics in Agriculture}, 175:105542, 2020.

\bibitem{quellec2020automatic}
Gwenol{\'e} Quellec, Mathieu Lamard, Pierre-Henri Conze, Pascale Massin, and
  B{\'e}atrice Cochener.
\newblock Automatic detection of rare pathologies in fundus photographs using
  few-shot learning.
\newblock {\em Medical image analysis}, 61:101660, 2020.

\bibitem{8935497}
Debasmit Das and C.~S.~George Lee.
\newblock A two-stage approach to few-shot learning for image recognition.
\newblock {\em IEEE Transactions on Image Processing}, 29:3336--3350, 2020.

\bibitem{9311786}
Xiaokang Zhou, Wei Liang, Shohei Shimizu, Jianhua Ma, and Qun Jin.
\newblock Siamese neural network based few-shot learning for anomaly detection
  in industrial cyber-physical systems.
\newblock {\em IEEE Transactions on Industrial Informatics}, 17(8):5790--5798,
  2021.

\bibitem{das2020few}
Debasmit Das, JH~Moon, and George Lee.
\newblock Few-shot image recognition with manifolds.
\newblock In {\em International Symposium on Visual Computing}, pages 3--14.
  Springer, 2020.

\bibitem{tibshirani1996regression}
Robert Tibshirani.
\newblock Regression shrinkage and selection via the lasso.
\newblock {\em Journal of the Royal Statistical Society: Series B
  (Methodological)}, 58(1):267--288, 1996.

\bibitem{markovsky2012low}
Ivan Markovsky.
\newblock {\em Low rank approximation: algorithms, implementation,
  applications}, volume 906.
\newblock Springer, 2012.

\bibitem{shetta2021convex}
Omar Shetta, Mahesan Niranjan, and Srinandan Dasmahapatra.
\newblock Convex multi-view clustering via robust low rank approximation with
  application to multi-omic data.
\newblock {\em IEEE/ACM transactions on computational biology and
  bioinformatics}, 2021.

\bibitem{markovsky2010approximate}
Ivan Markovsky and Mahesan Niranjan.
\newblock Approximate low-rank factorization with structured factors.
\newblock {\em Computational {S}tatistics \& {D}ata {A}nalysis},
  54(12):3411--3420, 2010.

\bibitem{papailiopoulos2013sparse}
Dimitris Papailiopoulos, Alexandros Dimakis, and Stavros Korokythakis.
\newblock Sparse pca through low-rank approximations.
\newblock In {\em International Conference on Machine Learning}, pages
  747--755. PMLR, 2013.

\bibitem{lee1999learning}
Daniel~D Lee and H~Sebastian Seung.
\newblock Learning the parts of objects by non-negative matrix factorization.
\newblock {\em Nature}, 401(6755):788--791, 1999.

\bibitem{li2017robust}
Zechao Li, Jinhui Tang, and Xiaofei He.
\newblock Robust structured nonnegative matrix factorization for image
  representation.
\newblock {\em IEEE transactions on neural networks and learning systems},
  29(5):1947--1960, 2017.

\bibitem{raghu2017svcca}
Maithra Raghu, Justin Gilmer, Jason Yosinski, and Jascha Sohl-Dickstein.
\newblock Svcca: Singular vector canonical correlation analysis for deep
  learning dynamics and interpretability.
\newblock {\em Advances in neural information processing systems}, 30, 2017.

\bibitem{wu2021generalized}
Aming Wu, Suqi Zhao, Cheng Deng, and Wei Liu.
\newblock Generalized and discriminative few-shot object detection via
  svd-dictionary enhancement.
\newblock {\em Advances in Neural Information Processing Systems}, 34, 2021.

\bibitem{9550610}
He~Zhang and Lili Liang.
\newblock Res-svdnet: A metric learning method for few-shot image
  classification.
\newblock In {\em 2021 40th Chinese Control Conference (CCC)}, pages
  7400--7405, 2021.

\bibitem{foley1975optimal}
Donald~H. Foley and John~W Sammon.
\newblock An optimal set of discriminant vectors.
\newblock {\em IEEE Transactions on computers}, 100(3):281--289, 1975.

\bibitem{leuschner2019supervised}
Johannes Leuschner, Maximilian Schmidt, Pascal Fernsel, Delf Lachmund, Tobias
  Boskamp, and Peter Maass.
\newblock Supervised non-negative matrix factorization methods for maldi
  imaging applications.
\newblock {\em Bioinformatics}, 35(11):1940--1947, 2019.

\bibitem{li2015robust}
Zechao Li, Jing Liu, Jinhui Tang, and Hanqing Lu.
\newblock Robust structured subspace learning for data representation.
\newblock {\em IEEE transactions on pattern analysis and machine intelligence},
  37(10):2085--2098, 2015.

\bibitem{tang2020cart}
Rong Tang and Xiaojun Zhang.
\newblock Cart decision tree combined with boruta feature selection for medical
  data classification.
\newblock In {\em 2020 5th IEEE International Conference on Big Data Analytics
  (ICBDA)}, pages 80--84. IEEE, 2020.

\bibitem{huang2015supervised}
Samuel~H Huang.
\newblock Supervised feature selection: A tutorial.
\newblock {\em Artif. Intell. Res.}, 4(2):22--37, 2015.

\bibitem{kursa2010feature}
Miron~B Kursa and Witold~R Rudnicki.
\newblock Feature selection with the boruta package.
\newblock {\em Journal of statistical software}, 36:1--13, 2010.

\bibitem{li2021semi}
Zechao Li and Jinhui Tang.
\newblock Semi-supervised local feature selection for data classification.
\newblock {\em Science China Information Sciences}, 64(9):1--12, 2021.

\bibitem{li2015unsupervised}
Zechao Li and Jinhui Tang.
\newblock Unsupervised feature selection via nonnegative spectral analysis and
  redundancy control.
\newblock {\em IEEE Transactions on Image Processing}, 24(12):5343--5355, 2015.

\bibitem{fisher1936use}
Ronald~A Fisher.
\newblock The use of multiple measurements in taxonomic problems.
\newblock {\em Annals of eugenics}, 7(2):179--188, 1936.

\bibitem{boyd2004convex}
Stephen Boyd, Stephen~P Boyd, and Lieven Vandenberghe.
\newblock {\em Convex Optimization}.
\newblock Cambridge university press, 2004.

\bibitem{chelali2009linear}
Fatma~Zohra Chelali, A~Djeradi, and R~Djeradi.
\newblock Linear discriminant analysis for face recognition.
\newblock In {\em 2009 International Conference on Multimedia Computing and
  Systems}, pages 1--10. IEEE, 2009.

\bibitem{NIPS2000_f9d11525}
Daniel Lee and H.~Sebastian Seung.
\newblock Algorithms for non-negative matrix factorization.
\newblock In T.~Leen, T.~Dietterich, and V.~Tresp, editors, {\em Advances in
  Neural Information Processing Systems}, volume~13. MIT Press, 2001.

\bibitem{dong2021transferred}
Aimei Dong, Zhigang Li, and Qiuyu Zheng.
\newblock Transferred subspace learning based on non-negative matrix
  factorization for eeg signal classification.
\newblock {\em Frontiers in Neuroscience}, 15, 2021.

\bibitem{chen2021deep}
Zhikui Chen, Shan Jin, Runze Liu, and Jianing Zhang.
\newblock A deep non-negative matrix factorization model for big data
  representation learning.
\newblock {\em Frontiers in Neurorobotics}, page~93, 2021.

\bibitem{Brunet4164}
Jean-Philippe Brunet, Pablo Tamayo, Todd~R. Golub, and Jill~P. Mesirov.
\newblock Metagenes and molecular pattern discovery using matrix factorization.
\newblock {\em Proceedings of the National Academy of Sciences},
  101(12):4164--4169, 2004.

\bibitem{pytorch}
Pytorch, forward and backward function hooks—pytorch documentation.

\bibitem{snell2017prototypical}
Jake Snell, Kevin Swersky, and Richard Zemel.
\newblock Prototypical networks for few-shot learning.
\newblock {\em Advances in Neural Information Processing Systems}, 30, 2017.

\bibitem{breastcancer}
Andrew Janowczyk and Anant Madabhushi.
\newblock Deep learning for digital pathology image analysis: A comprehensive
  tutorial with selected use cases.
\newblock {\em Journal of pathology informatics}, 7, 2016.

\bibitem{cruz2014automatic}
Angel Cruz-Roa, Ajay Basavanhally, Fabio Gonz{\'a}lez, Hannah Gilmore, Michael
  Feldman, Shridar Ganesan, Natalie Shih, John Tomaszewski, and Anant
  Madabhushi.
\newblock Automatic detection of invasive ductal carcinoma in whole slide
  images with convolutional neural networks.
\newblock In {\em Medical Imaging 2014: Digital Pathology}, volume 9041, page
  904103. SPIE, 2014.

\bibitem{cheng_2017}
Jun Cheng.
\newblock brain tumor dataset, Apr 2017.

\bibitem{cheng2015enhanced}
Jun Cheng, Wei Huang, Shuangliang Cao, Ru~Yang, Wei Yang, Zhaoqiang Yun,
  Zhijian Wang, and Qianjin Feng.
\newblock Enhanced performance of brain tumor classification via tumor region
  augmentation and partition.
\newblock {\em PloS one}, 10(10):e0140381, 2015.

\bibitem{he2020sample}
Xuehai He, Xingyi Yang, Shanghang Zhang, Jinyu Zhao, Yichen Zhang, Eric Xing,
  and Pengtao Xie.
\newblock Sample-efficient deep learning for covid-19 diagnosis based on ct
  scans.
\newblock {\em medrxiv}, 2020.

\bibitem{DeepDRi}
The 1st diabetic retinopathy – classification of fundus images according to
  the severity level of diabetic retinopathy.

\bibitem{acevedo2020dataset}
Andrea Acevedo, Anna Merino, Santiago Alf{\'e}rez, {\'A}ngel Molina, Laura
  Bold{\'u}, and Jos{\'e} Rodellar.
\newblock A dataset of microscopic peripheral blood cell images for development
  of automatic recognition systems.
\newblock {\em Data in Brief, ISSN: 23523409, Vol. 30,(2020)}, 2020.

\bibitem{al2020dataset}
Walid Al-Dhabyani, Mohammed Gomaa, Hussien Khaled, and Aly Fahmy.
\newblock Dataset of breast ultrasound images.
\newblock {\em Data in brief}, 28:104863, 2020.

\bibitem{moon2020computer}
Woo~Kyung Moon, Yan-Wei Lee, Hao-Hsiang Ke, Su~Hyun Lee, Chiun-Sheng Huang, and
  Ruey-Feng Chang.
\newblock Computer-aided diagnosis of breast ultrasound images using ensemble
  learning from convolutional neural networks.
\newblock {\em Computer methods and programs in biomedicine}, 190:105361, 2020.

\bibitem{tschandl2018ham10000}
Philipp Tschandl, Cliff Rosendahl, and Harald Kittler.
\newblock The ham10000 dataset, a large collection of multi-source
  dermatoscopic images of common pigmented skin lesions.
\newblock {\em Scientific data}, 5(1):1--9, 2018.

\bibitem{codella2019skin}
Noel Codella, Veronica Rotemberg, Philipp Tschandl, M~Emre Celebi, Stephen
  Dusza, David Gutman, Brian Helba, Aadi Kalloo, Konstantinos Liopyris, Michael
  Marchetti, et~al.
\newblock Skin lesion analysis toward melanoma detection 2018: A challenge
  hosted by the international skin imaging collaboration (isic).
\newblock {\em arXiv preprint arXiv:1902.03368}, 2019.

\bibitem{kermany2018identifying}
Daniel~S Kermany, Michael Goldbaum, Wenjia Cai, Carolina~CS Valentim, Huiying
  Liang, Sally~L Baxter, Alex McKeown, Ge~Yang, Xiaokang Wu, Fangbing Yan,
  et~al.
\newblock Identifying medical diagnoses and treatable diseases by image-based
  deep learning.
\newblock {\em Cell}, 172(5):1122--1131, 2018.

\bibitem{bilic2019liver}
Patrick Bilic, Patrick~Ferdinand Christ, Eugene Vorontsov, Grzegorz Chlebus,
  Hao Chen, Qi~Dou, Chi~Wing Fu, Xiao Han, Pheng-Ann Heng, J{\"u}rgen Hesser,
  et~al.
\newblock The liver tumor segmentation benchmark ({LiTS}).
\newblock {\em arXiv preprint arXiv:1901.04056}, 2019.

\bibitem{medmnistv2}
Jiancheng Yang, Rui Shi, Donglai Wei, Zequan Liu, Lin Zhao, Bilian Ke,
  Hanspeter Pfister, and Bingbing Ni.
\newblock Medmnist v2: A large-scale lightweight benchmark for 2d and 3d
  biomedical image classification.
\newblock {\em arXiv preprint arXiv:2110.14795}, 2021.

\bibitem{kather2019predicting}
Jakob~Nikolas Kather, Johannes Krisam, Pornpimol Charoentong, Tom Luedde,
  Esther Herpel, Cleo-Aron Weis, Timo Gaiser, Alexander Marx, Nektarios~A
  Valous, Dyke Ferber, et~al.
\newblock Predicting survival from colorectal cancer histology slides using
  deep learning: A retrospective multicenter study.
\newblock {\em PLoS medicine}, 16(1):e1002730, 2019.

\bibitem{woloshuk2021situ}
Andre Woloshuk, Suraj Khochare, Aljohara~F Almulhim, Andrew~T McNutt, Dawson
  Dean, Daria Barwinska, Michael~J Ferkowicz, Michael~T Eadon, Katherine~J
  Kelly, Kenneth~W Dunn, et~al.
\newblock In situ classification of cell types in human kidney tissue using 3d
  nuclear staining.
\newblock {\em Cytometry Part A}, 99(7):707--721, 2021.

\bibitem{squires2017rank}
Steven Squires, Adam Pr{\"u}gel-Bennett, and Mahesan Niranjan.
\newblock Rank selection in nonnegative matrix factorization using minimum
  description length.
\newblock {\em Neural computation}, 29(8):2164--2176, 2017.

\bibitem{balasubramanian2002isomap}
Mukund Balasubramanian and Eric~L Schwartz.
\newblock The isomap algorithm and topological stability.
\newblock {\em Science}, 295(5552):7--7, 2002.

\bibitem{lake2019omniglot}
Brenden~M Lake, Ruslan Salakhutdinov, and Joshua~B Tenenbaum.
\newblock The omniglot challenge: a 3-year progress report.
\newblock {\em Current Opinion in Behavioral Sciences}, 29:97--104, 2019.

\bibitem{kingma2014adam}
Diederik~P Kingma and Jimmy Ba.
\newblock Adam: A method for stochastic optimization.
\newblock {\em arXiv preprint arXiv:1412.6980}, 2014.

\end{thebibliography}


\end{document}